\documentclass[10pt]{amsart}
\usepackage{amsmath}
\usepackage{amsxtra}
\usepackage{amscd}
\usepackage{amsthm}
\usepackage{amsfonts}
\usepackage{amssymb}
\usepackage{eucal}

\newtheorem{cor}[subsection]{Corollary}
\newtheorem{lem}[subsection]{Lemma}
\newtheorem{sublem}[subsection]{Sublemma}

\newtheorem{prop}[subsection]{Proposition}

\newtheorem{lemconstr}[subsection]{Lemma-Construction}

\newtheorem{thm}[subsection]{Theorem}

\theoremstyle{definition}

\theoremstyle{remark}

\newcommand{\lemconstrref}[1]{Lemma-Construction~\ref{#1}}
\newcommand{\thmref}[1]{Theorem~\ref{#1}}
\newcommand{\secref}[1]{Sect.~\ref{#1}}
\newcommand{\lemref}[1]{Lemma~\ref{#1}}
\newcommand{\sublemref}[1]{Sublemma~\ref{#1}}
\newcommand{\propref}[1]{Proposition~\ref{#1}}
\newcommand{\corref}[1]{Corollary~\ref{#1}}

\newcommand{\nc}{\newcommand}
\nc{\renc}{\renewcommand}
\nc{\ssec}{\subsection}
\nc{\sssec}{\subsubsection}
\nc{\on}{\operatorname}

\nc\ol{\overline}
\nc\wt{\widetilde}
\nc\tboxtimes{\wt{\boxtimes}}
\nc{\alp}{\alpha}

\nc{\ZZ}{{\mathbb Z}}
\nc{\NN}{{\mathbb N}}
\nc{\OO}{{\mathbb O}}
\renc{\SS}{{\mathbb S}}
\nc{\DD}{{\mathbb D}}
\nc{\GG}{{\mathbb G}}

\nc{\Fq}{{\mathbb F}_q}
\nc{\Fqb}{\ol{{\mathbb F}_q}}
\nc{\Ql}{\ol{{\mathbb Q}_\ell}}
\nc{\id}{\text{id}}
\nc\X{\mathcal X}

\nc{\Hom}{\on{Hom}}
\nc{\Lie}{\on{Lie}}
\nc{\Loc}{\on{Loc}}
\nc{\Pic}{\on{Pic}}
\nc{\Bun}{\on{Bun}}
\nc{\IC}{\on{IC}}
\nc{\Aut}{\on{Aut}}
\nc{\rk}{\on{rk}}
\nc{\Sh}{\on{Sh}}
\nc{\Perv}{\on{Perv}}
\nc{\pos}{{\on{pos}}}
\nc{\Conv}{\on{Conv}}
\nc{\Sph}{\on{Sph}}
\nc{\Sym}{\on{Sym}}
\nc{\BunBb}{\overline{\Bun}_B}
\nc{\Buno}{\overset{o}{\Bun}}
\nc{\BunPb}{{\overline{\Bun}_P}}
\nc{\BunBM}{\overline{\Bun}_{B(M)}}
\nc{\BunPbw}{{\widetilde{\Bun}_P}}
\nc{\BunBP}{\widetilde{\Bun}_{B,P}}
\nc{\GUb}{\overline{G/U}}
\nc{\GUPb}{\overline{G/U(P)}}

\nc{\Hhom}{\underline{\on{Hom}}}
\nc\syminfty{\on{Sym}^{\infty}}
\nc\lal{\ol{\lambda}}
\nc\xl{\ol{x}}
\nc\thl{\ol{\theta}}
\nc\nul{\ol{\nu}}
\nc\mul{\ol{\mu}}
\nc\Sum\Sigma
\nc{\oX}{\overset{o}{X}{}}

\nc{\M}{{\mathcal M}}
\nc{\N}{{\mathcal N}}
\nc{\F}{{\mathcal F}}
\nc{\D}{{\mathcal D}}
\nc{\Q}{{\mathcal Q}}
\nc{\Y}{{\mathcal Y}}
\nc{\G}{{\mathcal G}}
\nc{\E}{{\mathcal E}}
\nc{\CalC}{{\mathcal C}}
\nc\Dh{\widehat{\D}}

\nc{\C}{{\mathcal C}}
\nc{\K}{{\mathcal K}}
\renewcommand{\H}{{\mathcal H}}

\nc{\T}{{\mathcal T}}
\nc{\V}{{\mathcal V}}
\renc{\P}{{\mathcal P}}
\nc{\A}{{\mathcal A}}
\nc{\B}{{\mathcal B}}
\nc{\U}{{\mathcal U}}

\nc{\Gr}{\on{Gr}}

\nc{\frn}{{\check{\mathfrak u}(P)}}
\nc{\p}{\overline{\mathfrak p}}
\nc{\q}{\overline{\mathfrak q}}
\nc\f{{\mathfrak f}}

\nc{\qo}{{\mathfrak q}}
\nc{\po}{{\mathfrak p}}
\nc{\s}{{\mathfrak s}}
\nc\w{\text{w}}

\nc\mathi\iota
\nc\Spec{\on{Spec}}
\nc\Mod{\on{Mod}}
\nc{\tw}{\widetilde{\mathfrak t}}
\nc{\pw}{\widetilde{\mathfrak p}}
\nc{\qw}{\widetilde{\mathfrak q}}
\nc{\jw}{\widetilde j}

\nc{\grb}{\overline{\Gr}}
\nc{\I}{\mathcal I}

\nc{\lambdach}{{\check\lambda}}
\nc{\Lambdach}{{\check\Lambda}{}}
\nc{\much}{{\check\mu}}
\nc{\omegach}{{\check\omega}}
\nc{\nuch}{{\check\nu}}
\nc{\etach}{{\check\eta}}
\nc{\alphach}{{\check\alpha}}
\nc{\betach}{{\check\beta}}
\nc{\rhoch}{{\check\rho}}
\nc{\ch}{{\check h}}

\nc{\Hb}{\overline{\H}}

\nc{\BA}{{\mathbb{A}}}
\nc{\BC}{{\mathbb{C}}}
\nc{\BG}{{\mathbb{G}}}
\nc{\BH}{{\mathbb{H}}}
\nc{\BM}{{\mathbb{M}}}
\nc{\BN}{{\mathbb{N}}}
\nc{\BP}{{\mathbb{P}}}
\nc{\BR}{{\mathbb{R}}}
\nc{\BZ}{{\mathbb{Z}}}
\nc{\BV}{{\mathbb{V}}}
\nc{\BW}{{\mathbb{W}}}
\nc{\BU}{{\mathbb{U}}}
\nc{\BX}{{\mathbb{X}}}
\nc{\BY}{{\mathbb{Y}}}
\nc{\BK}{{\mathbb{K}}}

\nc{\CC}{{\mathcal{C}}}
\nc{\CA}{{\mathcal{A}}}
\nc{\CB}{{\mathcal{B}}}

\nc{\CE}{{\mathcal{E}}}
\nc{\CF}{{\mathcal{F}}}
\nc{\CG}{{\mathcal{G}}}
\nc{\CL}{{\mathcal{L}}}
\nc{\CM}{{\mathcal{M}}}
\nc{\CN}{{\mathcal{N}}}
\nc{\CO}{{\mathcal{O}}}
\nc{\CP}{{\mathcal{P}}}
\nc{\CQ}{{\mathcal{Q}}}
\nc{\CR}{{\mathcal{R}}}
\nc{\CS}{{\mathcal{S}}}
\nc{\CT}{{\mathcal{T}}}
\nc{\CU}{{\mathcal{U}}}
\nc{\CV}{{\mathcal{V}}}
\nc{\CW}{{\mathcal{W}}}
\nc{\CY}{{\mathcal{Y}}}
\nc{\CZ}{{\mathcal{Z}}}

\nc{\cM}{{\check{\mathcal M}}{}}
\nc{\csM}{{\check{\mathcal A}}{}}
\nc{\oM}{{\overset{\circ}{\mathcal M}}{}}
\nc{\obM}{{\overset{\circ}{\mathbf M}}{}}
\nc{\oCA}{{\overset{\circ}{\mathcal A}}{}}
\nc{\obA}{{\overset{\circ}{\mathbf A}}{}}
\nc{\ooM}{{\overset{\circ}{M}}{}}
\nc{\osM}{{\overset{\circ}{\mathsf M}}{}}
\nc{\vM}{{\overset{\bullet}{\mathcal M}}{}}
\nc{\nM}{{\underset{\bullet}{\mathcal M}}{}}
\nc{\oD}{{\overset{\circ}{\mathcal D}}{}}
\nc{\obD}{{\overset{\circ}{\mathbf D}}{}}
\nc{\oA}{{\overset{\circ}{\mathbb A}}{}}
\nc{\op}{{\overset{\bullet}{\mathbf p}}{}}
\nc{\cp}{{\overset{\circ}{\mathbf p}}{}}
\nc{\oU}{{\overset{\bullet}{\mathcal U}}{}}
\nc{\oZ}{{\overset{\circ}{\mathcal Z}}{}}
\nc{\ofZ}{{\overset{\circ}{\mathfrak Z}}{}}
\nc{\oF}{{\overset{\circ}{\fF}}}

\nc{\fa}{{\mathfrak{a}}}
\nc{\fb}{{\mathfrak{b}}}
\nc{\fg}{{\mathfrak{g}}}
\nc{\fgl}{{\mathfrak{gl}}}
\nc{\fh}{{\mathfrak{h}}}
\nc{\fj}{{\mathfrak{j}}}
\nc{\fm}{{\mathfrak{m}}}
\nc{\fn}{{\mathfrak{n}}}
\nc{\fu}{{\mathfrak{u}}}
\nc{\fp}{{\mathfrak{p}}}
\nc{\fq}{{\mathfrak{q}}}
\nc{\fr}{{\mathfrak{r}}}
\nc{\fs}{{\mathfrak{s}}}
\nc{\fsl}{{\mathfrak{sl}}}
\nc{\hsl}{{\widehat{\mathfrak{sl}}}}
\nc{\hgl}{{\widehat{\mathfrak{gl}}}}
\nc{\hg}{{\widehat{\mathfrak{g}}}}
\nc{\chg}{{\widehat{\mathfrak{g}}}{}^\vee}
\nc{\hn}{{\widehat{\mathfrak{n}}}}
\nc{\chn}{{\widehat{\mathfrak{n}}}{}^\vee}

\nc{\fA}{{\mathfrak{A}}}
\nc{\fB}{{\mathfrak{B}}}
\nc{\fD}{{\mathfrak{D}}}
\nc{\fE}{{\mathfrak{E}}}
\nc{\fF}{{\mathfrak{F}}}
\nc{\fG}{{\mathfrak{G}}}
\nc{\fK}{{\mathfrak{K}}}
\nc{\fL}{{\mathfrak{L}}}
\nc{\fM}{{\mathfrak{M}}}
\nc{\fN}{{\mathfrak{N}}}
\nc{\fP}{{\mathfrak{P}}}
\nc{\fU}{{\mathfrak{U}}}
\nc{\fV}{{\mathfrak{V}}}
\nc{\fZ}{{\mathfrak{Z}}}

\nc{\bb}{{\mathbf{b}}}
\nc{\bc}{{\mathbf{c}}}
\nc{\be}{{\mathbf{e}}}
\nc{\bj}{{\mathbf{j}}}
\nc{\bn}{{\mathbf{n}}}
\nc{\bp}{{\mathbf{p}}}
\nc{\bq}{{\mathbf{q}}}
\nc{\bu}{{\mathbf{u}}}
\nc{\bv}{{\mathbf{v}}}
\nc{\bx}{{\mathbf{x}}}
\nc{\by}{{\mathbf{y}}}
\nc{\bw}{{\mathbf{w}}}
\nc{\bA}{{\mathbf{A}}}
\nc{\bB}{{\mathbf{B}}}
\nc{\bC}{{\mathbf{C}}}
\nc{\bD}{{\mathbf{D}}}
\nc{\bG}{{\mathbf{G}}}
\nc{\bH}{{\mathbf{H}}}
\nc{\bY}{{\mathbf{Y}}}
\nc{\bT}{{\mathbf{T}}}
\nc{\bU}{{\mathbf{U}}}
\nc{\bL}{{\mathbf{L}}}
\nc{\bI}{{\mathbf{I}}}

\nc{\bM}{{\mathbf{M}}}
\nc{\bN}{{\mathbf{N}}}
\nc{\bV}{{\mathbf{V}}}
\nc{\bW}{{\mathbf{W}}}
\nc{\bX}{{\mathbf{X}}}
\nc{\bK}{{\mathbf{K}}}
\nc{\bZ}{{\mathbf{Z}}}

\nc{\sA}{{\mathsf{A}}}
\nc{\sB}{{\mathsf{B}}}
\nc{\sC}{{\mathsf{C}}}
\nc{\sD}{{\mathsf{D}}}
\nc{\sF}{{\mathsf{F}}}
\nc{\sK}{{\mathsf{K}}}
\nc{\sM}{{\mathsf{M}}}
\nc{\sO}{{\mathsf{O}}}
\nc{\sH}{{\mathsf{H}}}
\nc{\sQ}{{\mathsf{Q}}}
\nc{\sP}{{\mathsf{P}}}
\nc{\sZ}{{\mathsf{Z}}}
\nc{\sG}{{\mathsf{G}}}
\nc{\sfp}{{\mathsf{p}}}
\nc{\sr}{{\mathsf{r}}}
\nc{\sfb}{{\mathsf{b}}}
\nc{\sfc}{{\mathsf{c}}}
\nc{\sd}{{\mathsf{d}}}

\nc{\tA}{{\widetilde{\mathbf{A}}}}
\nc{\tB}{{\widetilde{\mathcal{B}}}}
\nc{\tg}{{\widetilde{\mathfrak{g}}}}
\nc{\tG}{{\widetilde{G}}}
\nc{\TM}{{\widetilde{\mathbb{M}}}{}}
\nc{\tO}{{\widetilde{\mathsf{O}}}{}}
\nc{\tU}{{\widetilde{\mathfrak{U}}}{}}
\nc{\TZ}{{\tilde{Z}}}
\nc{\tx}{{\tilde{x}}}
\nc{\tbv}{{\tilde{\bv}}}
\nc{\tfP}{{\widetilde{\mathfrak{P}}}{}}
\nc{\tz}{{\tilde{\zeta}}}
\nc{\tmu}{{\tilde{\mu}}}

\nc{\urho}{\underline{\rho}}
\nc{\uB}{\underline{B}}
\nc{\uC}{{\underline{\mathbb{C}}}}
\nc{\ui}{\underline{i}}
\nc{\uj}{\underline{j}}
\nc{\ofP}{{\overline{\mathfrak{P}}}}
\nc{\oB}{{\overline{\mathcal{B}}}}
\nc{\og}{{\overline{\mathfrak{g}}}}
\nc{\oI}{{\overline{I}}}

\nc{\eps}{\varepsilon}
\nc{\hrho}{{\hat{\rho}}}

\nc{\one}{{\mathbf{1}}}
\nc{\two}{{\mathbf{t}}}

\nc{\Rep}{{\mathop{\operatorname{\rm Rep}}}}
\nc{\Tot}{{\mathop{\operatorname{\rm Tot}}}}
\nc{\Ker}{{\mathop{\operatorname{\rm Ker}}}}
\nc{\Hilb}{{\mathop{\operatorname{\rm Hilb}}}}
\nc{\End}{{\mathop{\operatorname{\rm End}}}}
\nc{\Ext}{{\mathop{\operatorname{\rm Ext}}}}
\nc{\CHom}{{\mathop{\operatorname{{\mathcal{H}}\it om}}}}
\nc{\GL}{{\mathop{\operatorname{\rm GL}}}}
\nc{\gr}{{\mathop{\operatorname{\rm gr}}}}
\nc{\Id}{{\mathop{\operatorname{\rm Id}}}}
\nc{\de}{{\mathop{\operatorname{\rm def}}}}
\nc{\length}{{\mathop{\operatorname{\rm length}}}}
\nc{\supp}{{\mathop{\operatorname{\rm supp}}}}

\nc{\Cliff}{{\mathsf{Cliff}}}
\nc{\Fl}{{\mathsf{Fl}}}
\nc{\Fib}{{\mathsf{Fib}}}
\nc{\Coh}{{\mathsf{Coh}}}
\nc{\FCoh}{{\mathsf{FCoh}}}

\nc{\reg}{{\text{\rm reg}}}

\nc{\cplus}{{\mathbf{C}_+}}
\nc{\cminus}{{\mathbf{C}_-}}
\nc{\cthree}{{\mathbf{C}_*}}
\nc{\Qbar}{{\bar{Q}}}

\nc{\bh}{{\bar{h}}}
\nc{\bOmega}{{\overline{\Omega}}}

\nc{\seq}[1]{\stackrel{#1}{\sim}}

\nc{\bSet}{{\mathbf Set}}
\nc{\BSet}{{\mathbb Set}}
\nc{\BVect}{{\mathbb Vect}}
\nc{\wh}{\widehat}

\begin{document}

\title[Algebraic groups over a 2-dim
local field]{Algebraic groups over a 2-dimensional local field: \\
some further constructions}

\author{Dennis Gaitsgory and David Kazhdan  }

\date{June 2004}

\dedicatory{Dedicated to A.~Joseph on his 60th birthday}

\address{D.K.: Einstein Institute of Mathematics, The Hebrew University
of Jerusalem, Givat Ram, Jerusalem, 91904, Israel.}
\address{D.G: Department of Mathematics, The University of Chicago,
5734 University Ave., Chicago IL 60637, USA.}
\email{kazhdan@math.huji.ac.il; gaitsgde@math.uchicago.edu}

\maketitle

\section*{Introduction}

\ssec{}

Let $\bK$ be a local field, $G$ a split reductive group over $\bK$, and $G((t))$
the corresponding loop group, regarded as a group-indscheme. In \cite{GK}
we suggested a categorical framework in which one can study representations
of the group $G((t))(\bK)=G\left(\bK((t))\right)$.

\medskip

The main point is that $\BG:=G((t))(\bK)$ admits no interesting representations
on vector spaces, and we have to consider pro-vector spaces instead.
In more detail, we regard $\BG$ as a group-like object in the category
$\BSet:=\on{Ind}(\on{Pro}(\on{Ind}(\on{Pro}(Set_0))))$, where $Set_0$ denotes the
category of finite sets. We observe that $\BSet$ has a natural pseudo-action
on the category $\BVect=\on{Pro}(Vect)$ of pro-vector spaces, and we define
the category $\on{Rep}(\BG)$ to consist of pairs $(\BV,\rho)$, where
$\BV\in \BVect$, and $\rho$ is an action map $\BG\times \BV\to \BV$ in
the sense of the above pseudo-action, satisfying the usual properties.

\medskip

In \cite{GK} several examples of objects of $\on{Rep}(\BG)$ were considered.
One such example is the principal series representation $\Pi$,
considered by M.~Kapranov in \cite{Ka}. Combining the results of {\it loc. cit.}
and the formalism of adjoint functors developed in \cite{GK} we showed that
the endomorphism algebra of $\Pi$ could be identified with the Cherednik
double affine Hecke algebra.

\medskip

Another example is the "left regular" representation, corresponding
to functions on $\BG$, with respect to the action of $\BG$ on itself
by left translations, denoted $M(\BG)$. The main feature of $M(\BG)$
is that the right action develops an anomaly: instead of the action
of $\BG$ we obtain an action of the Kac-Moody central extension
$\wh\BG_0$ of $\BG$ by means of the multiplicative group $\bG_m$, 
induced by the adjoint action of $G$ on its Lie algebra.

\ssec{}

In the present paper we continue the study of the category $\on{Rep}(\BG)$.
It is natural to subdivide the contents into three parts:

\medskip

In the first part, which consists of Sections \ref{sect distr} and \ref{right adjoints},
we prove some general results about representability of various covariant functors
on the category $\on{Rep}(\BG)$. These results are valid when $\BG$ is replaced
by an arbitrary group-like object on $\BSet$. We also introduce the pro-vector space
of distributions on an object of $\BSet$ with values in a pro-vector space; this
notion is used in order to construct actions on invariants and coinvariants of
representations of $\BG$.

\medskip

The second part occupies Sections \ref{funct of coinv},
\ref{sect semiinv}, and \secref{proof of symmetry}. We study
representations of a central extension $\wh\BG$ of $\BG$
by means of $\bG_m$ with a fixed central character $c:\bG_m\to \BC^*$;
the corresponding category is denoted $\on{Rep}_c(\wh\BG)$,
and $(\wh\BG',c')$ denotes the opposite extension with its central
character, cf. \cite{GK}, Sect. 5.9.

Our goal here is to study the functor of semi-invariants
$$\underset{\BG}{\overset{\frac{\infty}{2}}\otimes}:
\on{Rep}_c(\wh\BG)\times \on{Rep}_{c'}(\wh\BG')\to \BVect,$$
which couples the categories of representations at opposite levels.
The motivation for the existence of such functor
is provided by the semi-infinite cohomology functor on
the category of representations of a Kac-Moody Lie algebra.

The construction of $\underset{\BG}{\overset{\frac{\infty}{2}}\otimes}$
presented here follows the categorical interpretation of
semi-infinite cohomology, developed by L.~Positselsky (unpublished).

We use the functor of semi-invariants to prove the main result of
the present paper, \thmref{endomorphisms}. This theorem describes
for any quasi pro-unipotent subgroup $\BH$ of $\BG$ (cf.
\secref{sect infl}) the ring of endomorphisms of the functor
$\on{Coinv}_\BH:\on{Rep}(\BG)\to \BVect$, as the algebra of
endomorphisms of a certain object in the category of representations
of $\wh\BG_0$.

In particular, we obtain a functorial
interpretation of the double affine (Cherednik) algebra
in terms of the category $\on{Rep}(\BG)$, as the algebra
of endomorphisms of the functor of coinvariants with
respect to the maximal quasi pro-unipotent subgroup
of $\BG$.

\medskip

The third part consists of Sections \ref{sect thick} and \ref{sect
regular}, preceded by some preliminaries in \secref{stacks}. We
construct some more examples of objects of $\on{Rep}(\BG)$, this
time using the moduli stack of bundles on an algebraic curve $X$
over $\bK$, when we think of the variable $t$ as a local coordinate
near some point $\bx\in X$.

In particular, we show in \thmref{regular} that in this way one
naturally produces a pro-vector space, endowed with an action of
$\BG\times \BG$, such the space of bi-coinvariants with respect to
the maximal quasi pro-unipotent subgroup $\bI^{00}$ of $\BG$ is a
bi-module over Cherednik's algebra, isomorphic to the regular
representation of this algebra.

\ssec{Notation}

We keep the notations introduced in \cite{GK}. In particular,
for a category $\CC$ we denote by $\on{Ind}(\CC)$
(resp., $\on{Pro}(\CC)$) its ind- (resp., pro-) completion.

For a filtering set $I$ and a collection $A_i$ of objects of $\CC$
indexed by $I$, we will denote by
$\underset{\underset{I}\longrightarrow}{"lim"}\, A_i$ the resulting
object of $\on{Ind}(\CC)$ and by
$\underset{\underset{I}\longrightarrow}{lim}\, A_i:=limInd
(\underset{\underset{I}\longrightarrow}{"lim"}\, A_i)\in \CC$ the
inductive limit of the latter, if it exists. The notation for
inverse families is similar.

As was mentioned above $Set_0$ denotes the category of finite sets. We use the
short-hand notation $\bSet=\on{Ind}(\on{Pro}(Set_0))$ and $\BSet=\on{Ind}(\on{Pro}(\bSet))$.
We denote by $Vect_0$ the category of finite-dimensional vector space,
$Vect\simeq \on{Ind}(Vect_0)$ is the category of vector spaces, and $\BVect:=\on{Pro}(Vect)$
is the category of pro-vector spaces.

\ssec{Acknowledgements}

We would like to thank A.~Shapira for pointing out two mistakes in the previous
version of the paper.

The research of D.G. is supported by a long-term fellowship at the
Clay Mathematics Institute and a grant from DARPA. He would also
like to thank the Einstein Institute of Mathematics of the Hebrew
University of Jerusalem and IHES, where this work was written. The
research of D.K. is supported by an ISF grant.

\ssec{A correction to \cite{GK}}   \label{correction}

As was pointed out by A.~Shapira, Lemma 2.13 of \cite{GK} is wrong.
Namely, he explained to us a counter-example of a pro-vector space
$\BV$, acted on by a discrete set $X$ (thought of as an object of $\bSet$),
such that the action of every element of $X$ on $\BV$ is trivial, whereas
the action of $X$ on $\BV$ in the sense of the pseudo-action of
$\bSet\subset \BSet$ on $\BVect$ is non-trivial. 
Namely, $\BV=\underset{n\in \BN}{\underset{\longleftarrow}{"lim"}}\, \on{Funct}_c(\BZ^{\geq n})$
and $X=\BN$, such that $i\in \BN$ acts on each $\on{Funct}_c(\BZ^{\geq n})$ by
$$
\begin{cases}
& f(x_n,x_{n+1},...)\mapsto f(x_n,x_{n+1},...,x_i+1,...)-f(x_n,x_{n+1},...,x_i,...), \,\, i\geq n \\
& f(x_n,x_{n+1},...)\mapsto 0, \,\, i< n.
\end{cases}
$$

However, we have the following assertion. Let $\BG$ be as in \cite{GK}, Sect. 1.12
let and $\Pi_1=(\BV_1,\rho_1)$, $\Pi_2=(\BV_2,\rho_2)$ be two objects of
$\Rep(\BG,\BVect)$. Assume that $\BV_1$ is strict as a pro-vector space,
i.e., that it can be represented as 
$\underset{\longleftarrow}{"lim"}\, \bV_1^i$, where the maps in the inverse systesm
$\bV_1^j\to \bV_1^i$ are surjective. Let $\phi:\BV_1\to \BV_2$ be a map in $\BVect$,
which intertwines the actions of the set $G({\mathbf F})=\BG^{top}$ on $\BV_1$ and $\BV_2$.

\begin{lem}   \label{2.13}
Under the above circumstances, the map $\phi$ is a map in $\Rep(\BG,\BVect)$.
\end{lem}

\begin{proof}

We will prove a more general assertion, when we do not require $\BV_1$ and $\BV_2$
to be representations of $\BG$ on $\BVect$, but just objects of endowed with an
action of $\BG$, regarded as an object of $\BSet$. We claim that a map $\BV_1\to \BV_2$
compatible with a point-wise action of $\BG^{top}$ is compatible with an action 
of $\BG$ as an object of $\BSet$, under the assumption that $\BV_1$ is strict.

We represent $\BG$ as $\underset{\longrightarrow}{"lim"}\, \BX_k$, $\BX_k\in \on{Pro}(\bSet)$,
and for each $k$, $\BX_k\simeq \underset{\longleftarrow}{"lim"}\, \bX^l_k$,
such the maps $(\bX^{l'}_k)^{top}\to (\bX^l_k)^{top}$ are surjective. 
The assertion of the lemma reduces immediately to the case when $\BV_2=\bW\in Vect$,
and $\BG$ is replaced by $\BX_k$. In this case
$$\CHom(\BX_k\otimes \BV_2,\bW)\simeq
\underset{i}{\underset{\longrightarrow}{lim}}\, \CHom(\BX_k\otimes \bV_2^i,\bW).$$
However, by the assumption on the inverse system $\{\bV_1^i\}$, for every $i$ the map
$$\Hom((\BX_k)^{top} \times \bV_1^i,\bW)\to \Hom((\BX_k)^{top} \times \BV_1,\bW)$$
is injective. This reduces us to the case when $\BV_1=\bV$ is an object of $Vect$. The rest of the proof
proceeds as in Lemma 2.13 of \cite{GK}.

\end{proof}

\section{The pro-vector space of distributions} \label{sect distr}

\ssec{}

Let $\BX$ be an object of $\BSet$ and $\BV\in \BVect$. Consider
the covariant functor on $\BVect$ that assigns to $\BW$ the set
of actions $\BX\times \BV\to \BW$. We claim that this functor
is representable. We will denote the representing object by
${\mathbb Distr}_c(\BX,\BV)\in \BVect$; its explicit construction is given below.
It is clear from the definition that covariant functor
$\BV \to {\mathbb Distr}_c(\BX,\BV)$ is right exact.

\medskip

We begin with some preliminaries of categorical nature:

\begin{lem}  \label{inductive limits}
The category $\BVect$ is closed under inductive limits.
\end{lem}

\begin{proof}

Since $\BVect$ is abelian, it is enough to show that it
is closed under direct sums.

Let $\BV^\kappa$ be a collection of pro-vector spaces,
$\BV^\kappa\simeq \underset{\longleftarrow}{"lim"}\, 
\bV^\kappa_{i^\kappa}$ with $i^\kappa$
running over a filtering set $I^\kappa$. Consider the set 
$\underset{\kappa}\Pi\, I^\kappa$,
whose elements can be thought of as families 
$\{\varphi(\kappa)\in I^\kappa, \forall \kappa\}$.
This set is naturally filtering, and
$$\underset{\kappa}\oplus \BV^\kappa \simeq \underset{\longleftarrow}{"lim"}\,
\left(\underset{\kappa}\oplus \bV^\kappa_{\phi(\kappa)} \right),$$
where the inverse system is taken with respect to
$\underset{\kappa}\Pi\, I^\kappa$.

\end{proof}

\ssec{}  \label{various distr}

Let us now describe explicitly the pro-vector space ${\mathbb Distr}_c(\BX,\BV)$,

\medskip

If $X$ is a finite set and $V$ is a finite-dimensional vector space,
let $\on{Distr}_c(X,V)$ be the set of $V$-valued functions on $X$,
thought of as distributions. If $\bX^0\in \on{Pro}(Set_0)$ equals
$\underset{\longleftarrow}{"lim"}\, X_i$ with $X_i\in Set_0$ and $V$
is as above, set $${\mathbb
Distr}_c(\bX^0,V)=\underset{\longleftarrow}{"lim"}\,
\on{Distr}_c(X_i,V)\in \BVect.$$ Set also
$\on{Distr}(\bX^0,V)=\underset{\longleftarrow}{lim}\,
\on{Distr}_c(X_i,V)\in Vect$, i.e.
$$\on{Distr}(\bX^0,V)=limProj\, {\mathbb Distr}_c(\bX^0,V).$$

If $\bX$ is an object of $\bSet$ equal to $\underset{\longrightarrow}{"lim"}\, \bX^j$,
$\bX^j\in \on{Pro}(Set_0)$ and $\bV\in Vect$ is
$\underset{\longrightarrow}{"lim"}\, V_m$ with $V_m\in Vect_0$,
set
$${\mathbb Distr}_c(\bX,V)=
\underset{\underset{j,m}{\longrightarrow}}{lim}\, 
{\mathbb Distr}_c(\bX^j,V_m)\in \BVect,$$ where the
inductive limit is taken in $\BVect$. Set also
$$\on{Distr}_c(\bX,V)=
\underset{\underset{j,m}{\longrightarrow}}{lim}\, \on{Distr}_c(\bX^j,V_m)\in Vect.$$
When $V$ is finite-dimensional, the latter is the vector space, which is the topological 
dual of the topological vector
space $\on{Funct}^{lc}(\bX,V^*)$ of locally constant functions on $\bX$ with values
in $V^*$. Note that $\on{Distr}_c(\bX,V)$ is not isomorphic to
$limProj\, {\mathbb Distr}_c(\bX,V)$ even if $V$ is finite-dimensional.

\medskip

For $\BX^0\in \on{Pro}(\bSet)$ equal to $\underset{\longleftarrow}{"lim"}\, \bX_l$
with $\bX_l\in \bSet$ and $\BV$ is a pro-vector space equal to
$\underset{\longleftarrow}{"lim"}\, \bV_n$, set
$${\mathbb Distr}_c(\BX^0,\BV)=
\underset{\underset{l,n}{\longleftarrow}}{lim}\, {\mathbb Distr}_c(\bX_l,\bV_n)\in 
\BVect.$$
Finally, for $\BX\in \BSet$ equal to $\underset{\longrightarrow}{"lim"}\, 
\BX^k$ and $\BV\in \BVect$,
set
$${\mathbb Distr}_c(\BX,\BV)=\underset{\underset{k}{\longrightarrow}}{lim}\,
{\mathbb Distr}_c(\BX^k,\BV).$$

\begin{lemconstr}
For ${\mathbb Distr}_c(\BX,\BV)\in \BVect$ constructed above, there exists a natural
isomorphism
$$\Hom_{\BVect}({\mathbb Distr}_c(\BX,\BV),\BW)\simeq {\mathcal Hom}
(\BX\otimes \BV,\BW).$$
\end{lemconstr}

\begin{proof}

By the definition of both sides, we can assume that $\BX\in \on{Pro}(\bSet)$ and
$\BW=\bW\in Vect$. We have the following (evident) sublemma:

\begin{sublem} \label{hom and proj}
If $\BU=\underset{\longleftarrow}{lim}\, \BU_m$, where the projective limit
is taken in the category $\BVect$, then for any
$\BX\in \on{Pro}(\bSet)$ and $\bW\in Vect$,
$${\mathcal Hom}(\BX\otimes \BU,\bW)\simeq
\underset{\longrightarrow}{lim}\, {\mathcal Hom}(\BX\otimes
\BU_m,\bW).$$
\end{sublem}

The sublemma implies that  we can assume that $\BV=\bV\in Vect$. By applying again the
construction of ${\mathbb Distr}_c(\BX,\bV)$, 
we reduce the assertion of the lemma further
to the case when $\BX=\bX\in \bSet$, i.e., we have to show that
$$\Hom_{\BVect}({\mathbb Distr}_c(\bX,\bV),\bW)\simeq
{\mathcal Hom}(\bX\otimes \bV,\bW).$$

By the construction of ${\mathbb Distr}_c(\bX,\bV)$ and the
definition of the action, we can assume that $\bX\in
\on{Pro}(Set_0)$ and $\bV$ is finite-dimensional. In this case the
assertion is evident.

\end{proof}

\noindent{\it Remark.} For fixed $\BX$ and $\BV$ as above we can
also consider the contravariant functor on $\BVect$, given by
$\BW\mapsto {\mathcal Hom}(\BX\times \BW,\BV)$. It is easy to see
that this functor is ind-representable, but \lemref{hom and ind lem}
shows that it is not in general representable. We will denote the
resulting object of $\on{Ind}(\BVect)$ by ${\mathbb
Funct}(\BX,\BV)$.

\ssec{}

Let now $\BX,\BY$ be two objects of $\BSet$. The associativity
constraint of the pseudo-action of $\BSet$ and $\BVect$ gives rise to a map
\begin{equation} \label{map of distributions}
{\mathbb Distr}_c(\BX\times \BY,\BV)\to 
{\mathbb Distr}_c(\BX,{\mathbb Distr}_c(\BY,\BV)).
\end{equation}

Let us now recall the following definition from \cite{GK}, Sect. 2.10:

\medskip

\noindent An object $\BX\in \BSet$ is said to satisfy condition
(**) if it can be represented as $\underset{\longrightarrow}{"lim"}\, \BX_k$
with each $\BX_k\in \on{Pro}(\bSet)$ being weakly strict. We remind 
(cf. \cite{GK}, Sect. 1.10) that an object $\BX'\in \on{Pro}(\bSet)$ 
is said to be weakly strict if it can be represented as
$\underset{\longrightarrow}{"lim"}\, \bX'_i$, $\bX'_i\in \bSet$, such
that the transition maps $\bX'_i\to \bX'_j$ are {\it weakly surjective};
in the case of interest when all $\bX'_i$'s are locally compact, the latter
condition means that the map of topological spaces 
$\bX'{}^{\on{top}}_i\to \bX'_j{}^{\on{top}}$ has dense image.

As was shown in \cite{GK}, Sect. 2.12, if $G$ is an algebraic group 
over $\bK$, then the corresponding object $\BG\in \BSet$ satisfies
condition (**).

\begin{prop}  \label{surj of distr}
If $\BX\in \BSet$ satisfies condition (**), then the 
map in \eqref{map of distributions} is surjective.
\footnote{We are grateful to Alon Shapira who discovered
an error in the previous version of the paper, where
the (**) assumption on $\BX$ was omited.}
\end{prop}

This map is not in general an isomorphism. To construct
a counter-example, it suffices to take $\BV=\BC$--the $1$-dimensional
vector space, and $\BY$ a discrete set $Y\in Set\simeq \on{Ind}(Set_0)$, 
regarded as an object of $\BSet$ by means of $Set_0\to \on{Pro}(\bSet)$.

\begin{proof}

We need to show that for a pro-vector space $\BW$, the map
\begin{equation} \label{map of homs}
{\mathcal Hom}(\BX\otimes {\mathbb Distr}_c(\BY,\BV),\BW)\to
{\mathcal Hom}((\BX\times \BY)\otimes \BV,\BW)
\end{equation}
is injective. We will repeatedly use the facts 
that the functor $limInd:\on{Ind}(Vect)\to Vect$
is exact and the functor $limProj:\on{Pro}(Vect)\to Vect$ is left-exact.

By assumption, $\BX$ can be written as 
$\underset{\longrightarrow}{"lim"}\, \BX_k$ with $\BX_k\in \on{Pro}(\bSet)$
being weakly strict. Set also
$\BW=\underset{\longleftarrow}{"lim"}\, \bW_j$, $\bW_j\in Vect$. Both sides
of \eqref{map of homs} are projective limits over $k$ and $j$ of the
corresponding objects with $\BX$ replaced by $\BX_k$ and $\BW$ replaced by
$\bW_j$. So, we can assume that $\BX$ is a weakly strict object of
$\on{Pro}(\bSet)$ and $\BW=\bW\in Vect$.

Let us write now $\BY=\underset{\longrightarrow}{"lim"}\, \BY_{k'}$ with
$\BY_{k'}\in \on{Pro}(\bSet)$,
in which case ${\mathbb Distr}_c(\BY,\BV)\simeq
\underset{\longrightarrow}{lim}\, {\mathbb Distr}_c(\BY_{k'},\BV)$, and
$${\mathcal Hom}((\BX\times \BY)\otimes \BV,\BW)\simeq
\underset{\longleftarrow}{lim}\, {\mathcal Hom}((\BX\times
\BY_{k'})\otimes \BV,\BW).$$

\begin{lem}  \label{hom and ind lem}
If $\BU=\underset{\longrightarrow}{lim}\, \BU_m$, the inductive limit taking place
in $\BVect$, then for an object
$\BX\in \BSet$, satisfying condition (**), and $\BW\in \BVect$, the natural map
$${\mathcal Hom}(\BX\otimes \BU,\BW)\to
\underset{\longleftarrow}{lim}\, {\mathcal Hom}(\BX\otimes \BU_m,\BW)$$
is injective. If $\BX\in \bSet$, then this map is an isomorphism.
\end{lem}

\begin{proof}

As above, we can assume that $\BW=\bW\in Vect$, and $\BX$ is a weakly strict
object of $\on{Pro}(\bSet)$. Assume first that $\BX=\bX\in \bSet$.  
In this case the assertion of the lemma follows
from the description of inductive limits in $\BVect$ given in 
\lemref{inductive limits}.

Thus, let $\BX$ be represented as $\underset{\longleftarrow}{"lim"}\, \bX_l$,
$\bX_l\in \bSet$, with the transition maps $\bX_{l'}\to\bX_l$ being weakly 
surjective. Then
$${\mathcal Hom}(\BX\otimes \BU,\bW)\simeq 
\underset{l}{\underset{\longrightarrow}{lim}}\,
{\mathcal Hom}(\bX_l\otimes \BU,\bW)\simeq
\underset{l}{\underset{\longrightarrow}{lim}}\,
\underset{m}{\underset{\longleftarrow}{lim}}\, {\mathcal Hom}
(\bX_l\otimes \BU_m,\bW),$$ 
and
$$\underset{m}{\underset{\longleftarrow}{lim}}\, 
{\mathcal Hom}(\BX\otimes \BU_m,\bW)\simeq
\underset{m}{\underset{\longleftarrow}{lim}}\,
\underset{l}{\underset{\longrightarrow}{lim}}\, {\mathcal
Hom}(\bX_l\otimes \BU_m,\bW).$$

However, by the assumption, the transition maps 
${\mathcal Hom}(\bX_l\otimes \BU_m,\bW)\to 
{\mathcal Hom}(\bX_{l'}\otimes \BU_m,\bW)$ are injective. Therefore,
the natural map
$$\underset{l}{\underset{\longrightarrow}{lim}}\,
\underset{m}{\underset{\longleftarrow}{lim}}\, {\mathcal
Hom}(\bX_l\otimes \BU_m,\bW)\to
\underset{m}{\underset{\longleftarrow}{lim}}\,
\underset{l}{\underset{\longrightarrow}{lim}}\, {\mathcal
Hom}(\bX_l\otimes \BU_m,\bW)$$
is injective. 

\end{proof}

Hence, we are reduced to the case when $\BY$ is also an object of
$\on{Pro}(\bSet)$. Using \sublemref{hom and proj}, we reduce the
assertion further to the case when $\BV=\bV\in Vect$ and
$\BY=\bY\in \bSet$.

If $\BX=\underset{\longleftarrow}{"lim"}\, \bX_l$ then
both sides of \eqref{map of homs} are inductive limits over
$l$ of the corresponding objects with $\BX$ replaced by $\bX_l$.
Thus, from now on we will assume that $\BX=\bX\in \bSet$,
and we have  to show that the map
\begin{equation}  \label{map between simpler homs}
{\mathcal Hom}(\bX\otimes {\mathbb Distr}_c(\bY,\bV),\bW)\to
{\mathcal Hom}((\bX\times \bY)\otimes \bV,\bW)
\end{equation}
is injective,
where on the left-hand side ${\mathcal Hom}$ is 
understood in the sense of the pseudo-action
of $\bSet\subset \BSet$ on $\BVect$.

By applying \lemref{hom and ind lem}, we reduce the assertion to the case when 
$\bY\in \on{Pro}(Set_0)$ and $\bV$ is finite-dimensional. It is clear
that when $\bY$ belongs to $Set_0$, the map in
\eqref{map between simpler homs} is an isomorphism. Consider now the case
when $\bY=\underset{\longleftarrow}{"lim"}\, Y_i$ with
$Y_i\in Set_0$ and $\bX=\underset{\longrightarrow}{"lim"}\, \bX_n$ with
$\bX_n\in \on{Pro}(Set_0)$. Then, by \sublemref{hom and proj}
\begin{align*}
&{\mathcal Hom}(\bX\otimes {\mathbb Distr}_c(\bY,\bV),\bW)\simeq
\underset{i}{\underset{\longrightarrow}{lim}} \,{\mathcal Hom}
(\bX\otimes {\mathbb Distr}_c(Y_i,\bV),\bW) \simeq \\
&\underset{i}{\underset{\longrightarrow}{lim}} \,
\underset{n}{\underset{\longleftarrow}{lim}}\, {\mathcal
Hom}(\bX_n\otimes {\mathbb Distr}_c(Y_i,\bV),\bW)\simeq
\underset{i}{\underset{\longrightarrow}{lim}}
\,\underset{n}{\underset{\longleftarrow}{lim}}\, {\mathcal
Hom}((\bX_n\times Y_i)\otimes \bV,\bW).
\end{align*}

We also have and identification
$${\mathcal Hom}((\bX\times \bY)\otimes \bV,\bW)\simeq
\underset{n}{\underset{\longleftarrow}{lim}}\, 
{\mathcal Hom}((\bX_n\times \bY)\otimes \bV,\bW)\simeq
\underset{n}{\underset{\longleftarrow}{lim}}
\,\underset{i}{\underset{\longrightarrow}{lim}}\, 
{\mathcal Hom}((\bX_n\times Y_i)\otimes \bV,\bW).$$

Since $Y_i$ are finite sets, we can assume that the transition maps
$Y_{i'}\to Y_i$ are surjective. Therefore, the map
$$\underset{i}{\underset{\longrightarrow}{lim}}
\,\underset{n}{\underset{\longleftarrow}{lim}}\, {\mathcal
Hom}((\bX_n\times Y_i)\otimes \bV,\bW)\to
\underset{n}{\underset{\longleftarrow}{lim}}
\,\underset{i}{\underset{\longrightarrow}{lim}}\, {\mathcal
Hom}((\bX_n\times Y_i)\otimes \bV,\bW)$$
is injective.

\end{proof}

\ssec{}

As an application of \propref{surj of distr}, we will prove the following result.

Let $\rho:\BX\times \BV\to \BW$ be an action map. We can consider $\on{ker}(\rho)$
and $\on{coker}(\rho)$ as functors on $\BVect$:
$$\on{ker}(\rho)(\BU)=\{\phi:\BU\to \BV\,|\, \rho\circ\phi=0\} \text{ and }
\on{coker}(\rho)(\BU)=\{\psi:\BW\to \BU\,|\, \psi\circ \rho=0\}.$$

As in \cite{GK}, Proposition 2.8, one shows that $\on{coker}(\rho)$ is always
representable, and $\on{ker}(\rho)$ is representable if condition (**) 
is satisfied.

\begin{cor} \label{action on coinv}
Let $\BY\times \BV\to \BV$ and $\BY\times \BW\to \BW$ be actions
commuting in the natural sense with $\rho$. Then, if $\BY$ satisfies
(**), we have an action of $\BY$ on $\on{coker}(\rho)$, and if $\BX$
satisfies condition (**), we have an action of $\BY$ on $\on{ker}(\rho)$.
\end{cor}

This corollary will be used when $\BV=\BW$, and both $\BX=\BG$ and $\BY=\BH$ 
are group-like objects in $\BSet$, whose actions on $\BV$ commute. 
In this case we obtain that
$\BG$ acts on both invariants and coinvariants of $\BH$ on $\BV$.

\begin{proof}

Let us first prove the assertion about the cokernel.
Note that $\on{coker}(\rho)$ is isomorphic to the cokernel of the map
${\mathbb Distr}_c(\BX,\BV)\to \BW$ obtained from $\rho$. We need to show that
the composition
$${\mathbb Distr}_c(\BY,\BW)\to \BW\to \on{coker}(\rho)$$ factors through
${\mathbb Distr}_c(\BY,\on{coker}(\rho))$. By the right-exactness of
the functor ${\mathbb Distr}_c(\BY,\cdot)$,
$${\mathbb Distr}_c(\BY,\on{coker}(\rho))\simeq
\on{coker}\bigl({\mathbb Distr}_c(\BY,{\mathbb Distr}_c(\BX,\BV))\to
{\mathbb Distr}_c(\BY,\BW)\bigr),$$ 
and it is enough to show that the composition
$${\mathbb Distr}_c(\BY,{\mathbb Distr}_c(\BX,\BV))\to \BW\to \on{coker}(\rho)$$
vanishes.

However, using \propref{surj of distr}, we can replace
${\mathbb Distr}_c(\BY,{\mathbb Distr}_c(\BX,\BV))$ by
${\mathbb Distr}_c(\BY\times\BX,\BV)$, and the required assertion follows
from the commutative diagram:
$$
\CD
{\mathbb Distr}_c(\BX,\BV) @>\rho>>  \BW @>>> \on{coker}(\rho)  \\
@AAA  @AAA  \\
{\mathbb Distr}_c(\BY\times\BX,\BV) @>\rho>> {\mathbb Distr}_c(\BY,\BW).
\endCD
$$

The proof for $\on{ker}(\rho)$ is similar. We have to show that the composition
$${\mathbb Distr}_c(\BX,{\mathbb Distr}_c(\BY,\on{ker}(\rho)))\to
{\mathbb Distr}_c(\BX,\BV)\to \BW$$
vanishes. Using \propref{surj of distr}, it is sufficient to show that the
composition
$${\mathbb Distr}_c(\BX\times \BY,\on{ker}(\rho))\to
{\mathbb Distr}_c(\BX,\BV)\to \BW$$
vanishes, which follows from the assumption.

\end{proof}

\section{Existence of certain left adjoint functors}  \label{right adjoints}

\ssec{}

In what follows $\BG$ will be group-like object in $\BSet$ satisfying assumption (**).
Following \cite{GK}, we will denote by $\on{Rep}(\BG)$ the category of representations 
of $\BG$ on $\BVect$.

\begin{prop}   \label{free}
The forgetful functor $\on{Rep}(\BG)\to \BVect$ admits a left adjoint.
\end{prop}

\begin{proof}

We have to prove for any $\BW\in \BVect$ the representability of the
functor on $\on{Rep}(\BG)$ given by $\Pi=(\BV,\rho)\mapsto
\on{Hom}_{\BVect}(\BW,\BV)$. This functor obviously commutes with
projective limits in $\on{Rep}(\BG)$; so, by Proposition 1.2 of
\cite{GK} (with Ind replaced by Pro), it is enough to show that it is 
pro-representable.

Consider the category of pairs $(\Pi,\alpha)$, where 
$\Pi=(\BV,\rho)$ is an object in $\on{Rep}(\BG)$ and $\alpha:\BW\to \BV$ is a map in $\BVect$.
For any such pair we obtain an action map $\BG\times \BW\to \BV$, and hence a map
${\mathbb Distr}_c(\BG,\BW)\to \BV$. Since for an object of $\BVect$ the class of its
quotient objects is clearly a set, the sub-class of those $(\Pi,\alpha)$, for which the above
map ${\mathbb Distr}_c(\BG,\BW)\to \BV$ is surjective, 
is also a set. This set is naturally
filtered, and let us denote it by $A(\BW)$; 
it is endowed with a functor to $\on{Rep}(\BG)$
given by $(\Pi,\alpha)\mapsto \Pi$.

We claim that $\underset{\underset{(\Pi,\alpha)\in A(\BW)}\longleftarrow}{lim}\, \Pi$
is the object on $\on{Pro}(\on{Rep})$, which pro-represents our functor.

Indeed, for $\Pi'=(\BV',\rho')\in \on{Rep}(\BG)$, the map
$$\on{Hom}_{\on{Pro}(\on{Rep}(\BG))}
(\underset{\underset{(\Pi,\alpha)\in A(\BW)}\longleftarrow}{lim}\, \Pi,\Pi')=
\underset{\underset{(\Pi,\alpha)\in A(\BW)}\longrightarrow}{lim}\,
\on{Hom}_{\on{Rep}(\BG)}(\Pi,\Pi')
\to \on{Hom}_{\BVect}(\BW,\BV')$$ is evident. Vice versa, 
given a map
$\BW\to \BV'$ consider the induced map ${\mathbb Distr}_c(\BG,\BW)\to \BV'$,
and let $\BU$ be its image.  We claim that the action map 
$\BG\times \BU\to \BV'$ factors
through $\BU$; this would mean that $\Pi:=(\BU,\rho'|_{\BU})$ 
is a sub-object of $\Pi'$,
and we obtain a morphism from $\underset{\underset{(\Pi,\alpha)
\in A(\BW)}\longleftarrow}{lim}\, 
\Pi$ to $\Pi'$.

Consider the commutative diagram:
$$
\CD
{\mathbb Distr}_c(\BG\times \BG,\BW) @>>> 
{\mathbb Distr}_c(\BG,\on{Distr}_c(\BG,\BW))  \\
@V{\on{mult}}VV   @VVV  \\
{\mathbb Distr}_c(\BG,\BW)    &   &{\mathbb Distr}_c(\BG,\BU)  \\
@VVV   @VVV  \\
\BU @>>> \BV'
\endCD
$$
We need to show that the image of the vertical map $\on{Distr}_c(\BG,\BU)\to \BV'$
is contained in $\BU$. Since, by construction, the morphism
${\mathbb Distr}_c(\BG,\BW)\to \BU$ is surjective, and the functor
${\mathbb Distr}_c(\BG,\cdot)$ is right-exact, it suffices to show that the image
of the composed vertical map is contained in $\BU$.

However, by  \propref{surj of distr}, it is sufficient to check that the composed map
$$\on{Distr}_c(\BG\times \BG,\BW)\to \BV'$$ has its image contained in $\BU$, but this
follows from the above diagram.

\end{proof}

\ssec{}

Let us now derive some corollaries of \propref{free}. We will denote the left adjoint
constructed above by $\BV\mapsto \on{Free}(\BV,\BG)$.

\begin{cor} \label{coind from sgr}
Let $\BG_1\to \BG_2$ be a homomorphism of group-objects
of $\BSet$. Then the natural forgetful functor $\on{Rep}(\BG_2)\to
\on{Rep}(\BG_1)$ admits a left adjoint.
\end{cor}

\begin{proof}

Let $\Pi_1$ be an object of $\on{Rep}(\BG_1)$. The functor on $\on{Rep}(\BG_2)$
given by $\Pi\mapsto \on{Hom}_{\BG_1}(\Pi_1,\Pi)$ commutes
with projective limits. Therefore, by Lemma 1.2 of \cite{GK} it suffices to
show that it is pro-representable.

Let $\BV_1$ be the pro-vector space underlying $\Pi_1$. We have an injection
$\on{Hom}_{\BG_1}(\Pi_1,\Pi)\hookrightarrow \on{Hom}_{\BVect}(\BV_1,\BV)$,
where $\BV$ is the pro-vector space underlying $\Pi$.

By \propref{free} we know that the functor $\Pi\mapsto \on{Hom}_{\BVect}(\BV_1,\BV)$
is representable. Therefore, the assertion of the proposition follows from
Proposition 1.4 of \cite{GK}.

\end{proof}

We will denote the resulting functor $\on{Rep}(\BG_1)\to
\on{Rep}(\BG_2)$ by $\Pi\mapsto \on{Coind}^{\BG_2}_{\BG_1}(\Pi)$
and call it the coinduction functor.

\begin{cor}  \label{ind lim in rep}
The category $\on{Rep}(\BG)$ is closed under inductive limits.
\end{cor}

\noindent{\it Remark.}
Note that if $\BG=\bG$ is a group-object in $\bSet$, then the proof of
\lemref{inductive limits} shows that the category $\on{Rep}(\bG,\BVect)$
is closed under inductive limits. Moreover, the forgetful functor
$\on{Rep}(\bG,\BVect)\to \BVect$ commutes with inductive limits.

For an arbitrary $\BG\in \BSet$, the latter fact is not true, and we need
to resort to \propref{free} even to show the existence of inductive limits.
We will always have a surjection from the inductive limit of underlying
pro-vector spaces to the pro-vector space, underlying the inductive limit.

\begin{proof}

Let $\Pi_i=(\BV_i,\rho_i)$ be a filtering family of objects
of $\on{Rep}(\BG)$. Consider the covariant functor $\sF$ on
$\on{Rep}(\BG)$ given by
$$\Pi\mapsto \underset{\longleftarrow}{lim}\, \on{Hom}_{\on{Rep}(\BG)}(\Pi_i,\Pi).$$
Consider also the functor $\sF'$ that sends $\Pi=(\BV,\rho)$ to
$\underset{\longleftarrow}{lim}\, \on{Hom}_{\BVect}(\BV_i,\BV)$.

By \propref{free} and  \lemref{inductive limits}, the functor $\sF'$ is
representable.  Hence, by Proposition 1.4 of \cite{GK}, we conclude
that $\sF$ is pro-representable. Since $\sF$ obviously commutes with
projective limits in $\on{Rep}(\BG)$, it is representable by Lemma 1.2
of \cite{GK}.

\end{proof}

\ssec{Inflation}   \label{sect infl}

Let us call a group-object $\bH$ of $\bSet$ quasi-unipotent if
it can be presented as $\underset{\longrightarrow}{"lim"}\, \bH_i$,
where $\bH_i$ are group-objects of $\on{Pro}(Set_0)$ and transition
maps being homomorphisms, cf. \cite{GK}.

Let us call a group-object $\BH\in\on{Pro}(\bSet)$ quasi pro-unipotent
if it can be presented as $\underset{\longleftarrow}{"lim"}\, \bH^l$,
where $\bH^l$ are quasi-unipotent group-objects of $\bSet$,
and the transition maps $\bH^{l'}\to \bH^{l}$ being weakly surjective
homomorphisms, cf. \cite{GK}, Sect. 1.10.

According to Lemma 2.7 of \cite{GK}, if $\BH$ is quasi pro-unipotent,
the functor of $\BH$-coinvariants
$$\on{Coinv}_\BH:\on{Rep}(\BH,\BVect)\to \BVect$$
is exact.

\begin{prop}  \label{inf}
If $\BH$ is quasi pro-unipotent,
the functor $\on{Coinv}_\BH$ admits a left adjoint.
\end{prop}

We will refer to the resulting adjoint functor as "inflation", and denote it
by $\BV\mapsto \on{Inf}^\BH(\BV)$.

\ssec{Proof of \propref{inf}}

Let us first take $\bH$ to be a quasi-unipotent
group-object of $\bSet$, isomorphic to
$\underset{\longrightarrow}{"lim"}\, \bH_i$, where $\bH_i$ are
group-objects in $\on{Pro}(Set_0)$.

Let us show that for a vector space $\bV$, the functor
$\on{Rep}(\bH,Vect)\to Vect$ given by $\Pi\mapsto \on{Hom}(\bV,\Pi_\bH)$
is pro-representable.

\medskip

For an index $i$, consider the object
$\on{Coind}^\bH_{\bH_i}(\bV)\in \on{Rep}(\bH,\BVect)$,
where $\bV$ is regarded  as a trivial representation of $\bH_i$, and
$\on{Coind}$ is as in \corref{coind from sgr}. Using
Proposition 2.4 of \cite{GK}, we obtain that
$\on{Coind}^\bH_{\bH_i}(\bV)$ is a well-defined object of
$\on{Pro}(\on{Rep}(\bH,Vect))$, which pro-represents the
functor $\Pi\mapsto \Pi^{\bH_i}$.

Note that if $\bH$ is locally compact, and $\bH_i\subset \bH$
is open, then $\on{Coind}^\bH_{\bH_i}(\bV)$ belongs in fact to
$\on{Rep}(\bH,Vect)$, and is isomorphic to the space of
compactly supported $\bV$-valued distributions on $\bH/\bH_i$,
i.e., to the ordinary compact induction.

\medskip

Since $\bH_i$ is compact, we have $\Pi^{\bH_i}\simeq \Pi_{\bH_i}$.
Therefore, for $j>i$ we have natural maps
$$\on{Coind}^\bH_{\bH_j}(\bV)\to \on{Coind}^\bH_{\bH_i}(\bV).$$
Therefore, we can consider the object
$$\underset{\longleftarrow}{lim}\, \on{Coind}^\bH_{\bH_j}(\bV )\in
\on{Pro}(\on{Rep}(\bH,Vect)),$$
where the projective limit is taken in the category
$\on{Pro}(\on{Rep}(\bH,Vect))$.

For $\Pi\in \on{Rep}(\bH,Vect)$ we have:
$$\on{Hom}(\underset{\longleftarrow}{lim}\, \on{Coind}^\bH_{\bH_i}(\bV ),\Pi)
\simeq
\underset{\longrightarrow}{lim} \on{Hom}(\bV,\Pi_{\bH_i}).$$

Since $\Pi_\bH\simeq \underset{\longrightarrow}{lim}\, \Pi_{\bH_i}$,
the RHS of the above expression is not in general
isomorphic to $\on{Hom}(\bV,\Pi_\bH)$, except when $\bV$ is
finite-dimensional. In the latter case we set
$\on{Inf}^\bH(\bV):=\underset{\longleftarrow}{lim}\, \on{Coind}^\bH_{\bH_i}(\bV)$.

For general $\bV$, isomorphic to
$\underset{\longrightarrow}{lim}\,\bV_k$ with $\bV_k\in Vect_0$, we set
$$\on{Inf}^\bH(\bV)=\underset{\longrightarrow}{lim}\, \on{Inf}^\bH(\bV_k),$$
where the inductive limit is taken in $\on{Pro}(\on{Rep}(\bH,Vect))$,
cf. \lemref{inductive limits}.

\medskip

Now, the existence (and construction) of the functor $\on{Inf}^\BH$ follows
from Proposition 2.4 of \cite{GK}. Namely, if $\BH=\underset{\longleftarrow}{"lim"}\, \bH^l$
with $\bH^l$ being group-objects in $\bSet$ as above,
and $\BV=\underset{\longleftarrow}{"lim"}\, \bV_m$, we set
$$\on{Inf}^\BH(\BV)=\underset{l,m}{\underset{\longleftarrow}{lim}}\,
\on{Inf}^{\bH^l}(\bV_m),$$
where the projective limit is taken in the category
$\on{Rep}(\BH,\BVect)\simeq \on{Pro}(\on{Rep}(\BH,Vect))$, and each
$\on{Inf}^{\bH^l}(\bV_m)$ is regarded as a representation of $\BH$ via $\BH\to \bH^l$.

\section{The functor of coinvariants}   \label{funct of coinv}

\ssec{}

From now on we will assume that the group-like object $\BG$ is
obtained from a split reductive group $G$ over $\bK$, as in
\cite{GK}, Sect. 2.12. More generally, we will consider a central
extension $\wh{G}$ of $G((t))$ as in Sect. 2.14 of \cite{GK}, and
denote by $\on{Rep}_c(\wh\BG)$ the category of representations of
$\wh\BG$ at level $c$.

\medskip

Let $\BH$ be a quasi pro-unipotent group-object in $\on{Pro}(\bSet)$.
Let $\BH\to \BG$ be a homomorphism, and we will assume that
we are given a splitting of the induced extension $\wh\BG|_\BH$.
In particular, we have the forgetful functor
$\on{Rep}_c(\wh\BG)\to \on{Rep}(\BH,\BVect)$.

Consider the functor
$$\on{Rep}_c(\wh\BG)\to \BVect,$$ given by
$\Pi\mapsto \on{Coinv}_\BH(\Pi)$. Let $E(\BG,\BH)_c$ denote the algebra
of endomorphisms of this functor.

\medskip

\noindent{\it Remark.}
One can regard $E(\BG,\BH)_c$ as an analogue of the Hecke algebra of a
locally compact subgroup with respect to an open compact subgroup.
Indeed, if $\bG$ is a locally compact group-like object in $\bSet$ and
$\bH\subset \bG$ is open and compact, the corresponding Hecke algebra,
which by definition is the algebra of $\bH$-bi-invariant compactly supported
functions on $\bG$, can be interpreted both, as the algebra of endomorphisms
of the representation $\on{Coind}^\bG_\bH(\BC)$, where $\BC$ is the
trivial representation, and as the algebra of endomorphisms of the functor
$\Pi\mapsto \on{Coinv}_\bH(\Pi):\on{Rep}(\bG,Vect)\to Vect$.

\ssec{}

Recall now the representation $M_c(\BG)$, introduced in Sect. 5.6 of
\cite{GK}. According to the main theorem of {\it loc.cit.}, the
structure of $\wh\BG$-representation on $M_c(\BG)$ extends naturally
to a structure of $\wh\BG\times \wh\BG'$-representation, where 
$\wh\BG'$ is the group-object of $\BSet$ corresponding to the central
extension $\wh{G}'$ of $G((t))$, the latter being the Baer sum of $\wh{G}$
and the canonical extension $\wh{G}_0$, corresponding to the adjoint
action of $G$ on its Lie algebra. The action of $\wh\BG'$ of $M_c(\BG)$
has central character $c'$, given by the formula in Sect. 5.9 of \cite{GK}.

In what follows we will call objects of $\on{Rep}_{c'}(\wh\BG')$ 
"representations at the opposite level" to that of $\on{Rep}_c(\wh\BG)$.
We will refer to the $\wh{\BG}'$-action on $M_c(\BG)$ as the 
"right action".

\medskip

Using \corref{action on coinv}, by 
taking $\BH$-coinvariants with respect to $\BH$ mapping
to $\wh\BG'$, we obtain an object of $\on{Rep}_c(\wh\BG)$ which we will
denote by $M_c(\BG,\BH)$. By construction, we have a natural map
$$E(\BG,\BH)_{c'}\to \on{End}_{\BVect}(M_c(\BG,\BH)).$$
However, since the $\wh\BG$ and $\wh\BG'$ actions on $M_c(\BG)$ commute,
from \lemref{2.13} we obtain that endomorphisms of 
$M_c(\BG,\BH)$, resulting from the above map, commute with the $\wh\BG$-action.

Hence, we obtain a map
\begin{equation}   \label{map in question}
E(\BG,\BH)_{c'}\to \on{End}_{\on{Rep}_c(\wh\BG)}(M_c(\BG,\BH)).
\end{equation}

We will prove the following theorem:
\begin{thm}   \label{endomorphisms}
The map in \eqref{map in question} is an isomorphism.
\end{thm}

\ssec{}

Let us consider a few examples. Suppose first that the group
$\BH$ is trivial. As a corollary of \thmref{endomorphisms}
we obtain:

\begin{thm} \label{Hecke}
The algebra $E(\BG)_c$ of endomorphisms of the forgetful functor
$\on{Rep}_c(\wh\BG)\to \BVect$ is isomorphic to the algebra of
endomorphisms of the object $M_c(\BG)\in \on{Rep}_{c'}(\wh\BG')$.
\end{thm}

Let now $\BH$ be a {\it thick}
subgroup of $\bG[[t]]$ (see \cite{GK}, Sect. 2.12). Note that in this
case, the object $M_c(\BG,\BH)$ is isomorphic to the induced representation
$i^{\wh{\BG}}_\BH(\BC)$ of \cite{GK}, Sect. 3.3, where $\BC$ is the trivial
1-dimensional representation of $\BH$.

In particular, let us take $\BH$ to be $\bI^{00}$, the subgroup of $\bI$ equal to the kernel of
the natural map $\bI\to \bT\to \Lambda$,
where  $I\subset G[[t]]$ is the Iwahori subgroup and
$\Lambda$ is the lattice of cocharacters of $T$, regarded as a
quotient of $\bT$ by its maximal compact subgroup.

The corresponding induced representation $i^{\wh\BG}_\BH(\BC)$
is isomorphic to Kapranov's representation, denoted in Sect. 4
of \cite{GK} by $\BV_c$. Assume now that $G$ is semi-simple and
simply-connected. In this case it follows from Corollary 4.4 of \cite{GK}
that the algebra $\on{End}(\BV_c)$ is isomorphic
to the Cherednik algebra $\overset{\cdot\cdot}{\sH}_{q,c'}$. 
From \thmref{endomorphisms} we obtain:

\begin{cor}
The Cherednik algebra $\overset{\cdot\cdot}{\sH}_{q,c'}$ is isomorphic
to the algebra of endomorphisms of the functor
$\Pi\to \on{Coinv}_{\bI^{00}}(\Pi):\on{Rep}_c(\wh\BG)\to \BVect$.
\end{cor}

\ssec{}

Note that by combining \propref{inf} and \corref{coind from sgr}, we obtain that
the above functor $\on{Coinv}_\BH:\on{Rep}_c(\wh\BG)\to \BVect$
admits a left adjoint:
$$\BV\mapsto \on{Coind}^{\wh\BG}_\BH(\on{Inf}^\BH(\BV)).$$

Of course, the algebra of endomorphisms
of this functor is isomorphic to $E(\BG,\BH)^o_c$.

\medskip

Consider now the functor $\on{Rep}_c(\wh\BG)\to Vect$ obtained by
composing $\on{Coinv}_\BH$ with the functor $limProj: \BVect\to Vect$.
Let $\ol{E}(\BG,\BH)_c$ be the algebra of endomorphisms of this latter
functor. We have a natural map $E(\BG,\BH)_c\to \ol{E}(\BG,\BH)_c$.

\begin{prop}

\smallskip

\noindent {\em (a)} The map $E(\BG,\BH)_c\to \ol{E}(\BG,\BH)_c$ is injective.

\smallskip

\noindent{\em (b)}  The algebra $\ol{E}(\BG,\BH)_c^o$ is isomorphic to
$\on{End}_{\on{Rep}_c(\wh\BG)}\left(\on{Coind}^{\wh\BG}_\BH(\on{Inf}^\BH(\BC))\right)$.
\end{prop}

We do not know under what conditions on $\BH$ one might expect that
the above map $E(\BG,\BH)_c\to \ol{E}(\BG,\BH)_c$ is an isomorphism.

\begin{proof}

To prove the first assertion of the proposition, note that by \thmref{endomorphisms},
the evaluation map $E(\BG,\BH)_c\to 
\on{End}_{\BVect}\left(\on{Coinv}_\BH(M_c(\BG))\right)$ is injective.

By construction, the pro-vector space $M_c(\BG)$
can be represented as a countable inverse limit with surjective restriction maps. Hence,
by Proposition 2.5 of \cite{GK}, $\on{Coinv}_\BH(M_c(\BG))\in \BVect$ will also
have this property. We have:

\begin{lem}
For any pro-vector space, which can be represented as a countable 
inverse limit with surjective restriction maps, the morphism
$limProj(\BV)\to\BV$ is surjective.
\end{lem}

This lemma implies that the map 
$\on{End}_{\BVect}(\BV)\to \on{End}_{Vect}(limProj(\BV))$
is injective.

\medskip

To prove the second assertion, we must analyze the endomorphism algebra of the functor
$Vect\to \on{Rep}_c(\wh\BG)$ given by
$$\bV\mapsto \on{Coind}^{\wh\BG}_\BH(\on{Inf}^\BH(\bV)).$$

However, as every left adjoint, this functor commutes with inductive limits. Therefore,
its enough to consider its restriction to the subcategory $Vect_0$. This implies
the proposition.

\end{proof}

\section{The functor of semi-invariants}   \label{sect semiinv}

\ssec{}

Our method of proof of \thmref{endomorphisms} in based on considering the functor
of $\BG$-semi-invariants
$$\underset{\BG}{\overset{\frac{\infty}{2}}{\otimes}}:\on{Rep}_{c'}(\wh\BG')\times
\on{Rep}_c(\wh\BG)\to \BVect,$$
where $c$ and $c'$ are opposite levels. The construction of this functor
mimics the construction of the semi-infinite cohomology functor for associative
algebras by L.~Positselsky, \cite{Pos}.

\medskip

For $\Pi_c\in \on{Rep}_c(\wh\BG)$, $\Pi_{c'}\in \on{Rep}_{c'}(\wh\BG')$ consider
the pro-vector spaces
$$\Pi_{c'}\otimes \Pi_c \text{ and } \Pi_{c'}\otimes M_c(\BG)\otimes \Pi_c.$$
We consider the former as acted on by the diagonal copy of $\bG[[t]]$, and the
latter by two mutually commuting copies of $\bG[[t]]$: one acts
diagonally on $\Pi_{c'}\otimes M_c(\BG)$ via the {\it left} $\wh\BG$-action
on $M_c(\BG)$; the other copy acts diagonally on $M_c(\BG)\otimes \Pi_c$
via the {\it right} action. Consider the object
$$\left(\Pi_{c'}\otimes M_c(\BG)\otimes \Pi_c\right)_{\bG[[t]]\times \bG[[t]]}.$$
We will construct two natural maps
\begin{equation} \label{two maps}
\left(\Pi_{c'}\otimes \Pi_c\right)_{\bG[[t]]}\rightrightarrows
\left(\Pi_{c'}\otimes M_c(\BG)\otimes \Pi_c\right)_{\bG[[t]]\times \bG[[t]]}.
\end{equation}

To construct the first map recall from Lemma 5.8 of \cite{GK} that
\begin{equation} \label{id of ten prod}
\left(M_c(\BG)\otimes \Pi_c\right)_{\bG[[t]]}\simeq
i^{\wh\BG}_{\bG[[t]]}\left(r^{\wh\BG}_{\bG[[t]]}(\Pi_c)\right).
\end{equation}
Since $\BG/\bG[[t]]$ is ind-compact, the functor
$i^{\wh\BG}_{\bG[[t]]}$ is isomorphic to the induction functor,
$\wt{i}^{\wh\BG}_{\bG[[t]]}$. Therefore, we obtain a morphism
of $\wh\BG$-representations
\begin{equation} \label{induction map}
\Pi_c\to i^{\wh\BG}_{\bG[[t]]}\left(r^{\wh\BG}_{\bG[[t]]}(\Pi_c)\right)\simeq
\left(M_c(\BG)\otimes \Pi_c\right)_{\bG[[t]]}
\end{equation}
by adjunction from the identity map 
$r^{\wh\BG}_{\bG[[t]]}(\Pi_c)\to r^{\wh\BG}_{\bG[[t]]}(\Pi_c)$.

The first map in \eqref{two maps} comes from \eqref{induction map}
by tensoring with $\Pi_{c'}$ and taking $\bG[[t]]$-coinvariants.

\medskip

To construct the second map in \eqref{two maps} we will use the following
observation. Let $\wt{M}_c(\BG)$ be a representation of
$\wh\BG\times \wh\BG'$, obtained from the representation $M_{c'}(\BG)$
of $\wh\BG'\times \wh\BG$, by flipping the roles of $\wh\BG$ and $\wh\BG'$.
We have:

\begin{prop} \label{symmetry} \hfill

\smallskip

\noindent {\em (1)} We have a natural isomorphism 
of $\wh{\BG}\times \wh{\BG'}$-representations
$\wt{M}_c(\BG)\simeq M_{c'}(\BG)$.

\smallskip

\noindent {\em (2)}
The resulting two morphisms $$M_c(\BG)\rightrightarrows
\bigl(M_c(\BG)\otimes M_c(\BG)\bigr)_{\bG[[t]]}$$
one, coming from \eqref{induction map}, and the other from interchanging
the roles of $c$ and $c'$, coincide.
\end{prop}

\noindent{\it Remark.} It will follow from the proof, that statement (2)
of the proposition fixes the isomorphism of statement (1) uniquely.

\medskip

The proof will be given in \secref{proof of symmetry}. Using this
proposition we construct the second map in \eqref{two maps} by
simply interchanging the roles of $c$ and $c'$.

\ssec{}   \label{def semiinv}

For $\Pi_c,\Pi_{c'}$ as above, we set
$\Pi_{c'}\underset{\BG}{\overset{\frac{\infty}{2}}{\otimes}}\Pi_c$ to
be the equalizer (i.e., the kernel of the difference) of the two maps in \eqref{two maps}.
Note that since the functor of $\bG[[t]]$-coinvariants is only right-exact,
the resulting functor $\underset{\BG}{\overset{\frac{\infty}{2}}{\otimes}}$
is a priori neither right nor left exact.

Suppose now that $\Pi_c$ is not only a representation of $\wh\BG$,
but carries an additional commuting action of some group-object $\BH\in \BSet$,
which satisfies condition (**). In this case it follows from 
\corref{action on coinv}
that $\Pi_{c'}\underset{\BG}{\overset{\frac{\infty}{2}}{\otimes}}\Pi_c$
is an object of $\on{Rep}(\BH)$.

\medskip

The key assertion describing the behavior 
of the functor of semi-invariants is the following:

\begin{prop}   \label{computation}
For $M_c(\BG)$, regarded as an object of $\on{Rep}_c(\wh\BG)$, we have
a natural isomorphism
$\Pi_{c'}\underset{\BG}{\overset{\frac{\infty}{2}}{\otimes}}M_c(\BG)\simeq 
\Pi_{c'}$. Moreover,
this isomorphism is compatible with the $\wh\BG$-actions.
\end{prop}

\begin{proof}

Consider the following general set-up.
Let $\CC_1$ and $\CC_2$ be two abelian categories,
$\sG:\CC_1\to\CC_2$ be a functor, and $\sF:\CC_2\to \CC_1$
its right adjoint. By composing with $\sF\circ \sG$ on the left and on the right,
the adjunction map
$\on{Id}_{\CC_1}\to \sF\circ \sG$ gives rise to two maps
\begin{equation} \label{complex}
\sF\circ \sG \rightrightarrows \sF\circ \sG\circ \sF\circ \sG,
\end{equation}
such that $\on{Id}_{\CC_1}$  maps to thei equalizer.

\begin{lem}
Assume that the functor $\sG$ is exact and faithfull. Then the map
$$\on{Id}_{\CC_1}\to \on{Equalizer}\left(\sF\circ \sG \rightrightarrows \sF\circ \sG\circ \sF\circ \sG\right)$$
is an isomorphism.
\end{lem}

\begin{proof}
By assumption o $\sG$, it is enough to show that 
$$\sG\to \on{Equalizer}\left(\sG\circ \sF\circ \sG\rightrightarrows \sG\circ \sF\circ \sG\circ \sF\circ \sG\right)$$
is an isomorphism, but this happens for any pair of adjoint functors.
\end{proof}

We apply this lemma to $\CC_1=\on{Rep}_{c'}(\wh\BG')$,
$\CC_2=\on{Rep}(\bG[[t]],\BVect)$ with $\sF=i^{\wh\BG}_{\bG[[t]]}$,
$\sG=r^{\wh\BG}_{\bG[[t]]}$.
To prove the Proposition it is sufficient to show that for 
$\Pi_{c'}\in \on{Rep}_{c'}(\wh\BG')$ the terms and maps
in \eqref{two maps} are equal to the corresponding ones in
\eqref{complex}.

First, by \eqref{id of ten prod} and \propref{symmetry}(1), for $\Pi_{c'}$ as above,
$\sF\circ \sG(\Pi_{c'})$ is indeed isomorphic to $(\Pi_{c'}\otimes M_c(\BG))_{\bG[[t]]}$.
Furthermore, by applying the functor $\sF\circ \sG$ to the adjunction map
$\Pi_{c'}\to \sF\circ \sG(\Pi_{c'})$ we obtain the second of the two maps from
\eqref{two maps}.

Let us now calculate the adjunction map $\on{Id}_{\Rep_{c'}(\wh\BG')}\to
i^{\wh\BG}_{\bG[[t]]}\circ r^{\wh\BG}_{\bG[[t]]}$ applied to 
$$\sF\circ \sG(\Pi_{c'})\simeq (\Pi_{c'}\otimes M_c(\BG))_{\bG[[t]]}.$$ By construction,
it is obtained from the adjunction map
$$M_c(\BG)\to \sF\circ \sG(M_c(\BG))\simeq (M_c(\BG)\otimes M_c(\BG))_{\bG[[t]]}$$
by tensoring with $\Pi_{c'}$ and taking $G[[t]]$-coinvariants.
Therefore, by  \propref{symmetry}(2), it coincides with 
the first map from \eqref{two maps}.

\end{proof}

\noindent{\it Remark.} Note that by \propref{symmetry}(2), the two identifications
$M_c(\BG)\underset{\BG}{\overset{\frac{\infty}{2}}{\otimes}}M_c(\BG)\simeq M_c(\BG)$,
one coming from \propref{computation} applied to $\Pi_{c'}=M_c(\BG)$, and the other
from interchanging the roles of $c$ and $c'$ as in \propref{symmetry}(1), coincide.

\ssec{Proof of \thmref{endomorphisms}}  \label{proof of end}

Let $\Pi_{c'}$ be an object of $\on{Rep}_{c'}(\wh\BG')$, and let $\Pi_c$ be an object of
$\on{Rep}_c(\wh\BG)$, carrying an additional commuting action of a group-object $\BH\in\BSet$,
which is quasi pro-unipotent. Then, using \corref{action on coinv} and the fact that the functor
$\on{Coinv}_\BH$ is exact (Lemma 2.7 of \cite{GK}), we obtain an isomorphism:
$$\Pi_{c'}\underset{\BG}{\overset{\frac{\infty}{2}}{\otimes}} (\Pi_c)_{\BH}\simeq
\left(\Pi_{c'}\underset{\BG}{\overset{\frac{\infty}{2}}{\otimes}} \Pi_c\right )_{\BH}.$$

Applying this for $\Pi_c=M_c(\BG)$, we obtain a functorial isomorphism:
\begin{equation} \label{semiinv against ind}
\Pi_{c'}\underset{\BG}{\overset{\frac{\infty}{2}}{\otimes}}
M_c(\BG,\BH)\simeq (\Pi_{c'})_\BH.
\end{equation}

Therefore, we obtain a map
\begin{equation}  \label{other map}
\on{End}_{\on{Rep}_c(\wh\BG)}(M_c(\BG,\BH))\to E(\BG,\BH)_{c'}.
\end{equation}

The fact that the composition
$$\on{End}_{\on{Rep}_c(\wh\BG)}(M_c(\BG,\BH))\to E(\BG,\BH)_{c'}\to
\on{End}_{\on{Rep}_c(\wh\BG)}(M_c(\BG,\BH))$$ is the
identity map follows from the remark following the proof of \propref{computation}.

\medskip

Therefore, to finish the proof of the theorem it suffices to show that the map
of \eqref{map in question} is injective. For that note, that for any
$\Pi_{c'}\in \on{Rep}_{c'}(\wh\BG')$
we have an injection
$\Pi_{c'}\hookrightarrow (\Pi_{c'}\otimes M_c(\BG))_{\bG[[t]]}$
(coming from the above adjunction $\on{Id}_{\Rep_{c'}(\wh\BG')}\to
i^{\wh\BG}_{\bG[[t]]}\circ r^{\wh\BG}_{\bG[[t]]}$)
and a surjection
$\Pi_{c'}\otimes M_c(\BG)\twoheadrightarrow (\Pi_{c'}\otimes M_c(\BG))_{\bG[[t]]}$
of objects of $\on{Rep}_{c'}(\wh\BG')$.

\begin{lem}
Suppose an element $\alpha\in E(\BG,\BH)_{c'}$ annihilates
$(\Pi_{c'})_\BH$ for some $\Pi_{c'}\in \on{Rep}_{c'}(\wh\BG')$. Then
$\alpha$ annihilates all objects of the form $(\BV\otimes
\Pi_{c'})_\BH$ for $\BV\in \BVect$.
\end{lem}

\begin{proof}

Suppose that $\BV=\underset{\longleftarrow}{"lim"}\,\bV_i$,
$\bV_i\in Vect$. Then $\BV\otimes \Pi_{c'}\simeq
\underset{\longleftarrow}{lim}\, (\bV_i\otimes \Pi_{c'})$,
where the projective limit is taken in the category $\BVect$.

Using Corollary 2.6 of \cite{GK}, we have:
$(\BV\otimes \Pi_{c'})_\BH\simeq
\underset{\longleftarrow}{lim}\, (\bV_i\otimes \Pi_{c'})_\BH$.
This shows that we can assume that $\BV$ is a {\it vector space},
which we will denote by $\bV$.

Let us write $\bV=\underset{\longrightarrow}{lim}\, \bV_i$, where
$\bV_i\in Vect_0$.

\begin{sublem}
For $\bV=\underset{\longrightarrow}{lim}\, \bV_i$ and $\BW\in \BVect$
the natural map
$$\underset{\longrightarrow}{lim}\, (\bV_i\otimes \BW)\to
(\underset{\longrightarrow}{lim}\, \bV_i)\otimes \BW$$
is surjective.
\end{sublem}

Therefore, we have a surjection
$$\underset{\longrightarrow}{lim}\, (\bV_i\otimes \Pi_{c'})\twoheadrightarrow
\bV\otimes \Pi_{c'},$$ and, hence, a surjection on the level of
coinvariants.

Since by assumption, $\alpha$ annihilates every $(\bV_i\otimes \Pi_{c'})_\BH$,
and the functor $\on{Coinv}_\BH$ commutes with inductive limits
(cf. \corref{action on coinv}), we obtain that $\alpha$ annihilates also
$\left(\underset{\longrightarrow}{lim}\, (\bV_i\otimes \Pi_{c'})\right)_\BH$. Hence,
by the above, it annihilates also $(\bV\otimes \Pi_{c'})_\BH$.

\end{proof}

Using this lemma and the exactness of the functor of $\BH$-coinvariants,
we obtain that any $\alpha\in \on{ker}(E(\BG,\BH)_{c'}\to \on{End}(M_c(\BG,\BH))$
annihilates all $(\Pi_{c'}\otimes M_c(\BG))_\BH$, and hence $(\Pi_{c'})_\BH$
for any $\Pi_{c'}$.

\medskip

\noindent{\it Remark.} Note that the same argument proves the following more general 
assertion. Let $\BH_1,\BH_2$ be two quasi pro-unipotent groups endowed with
homomorphisms to $\wh\BG$.
Then the space of natural transformations between the functors
$\on{Coinv}_{\BH_1},\on{Coinv}_{\BH_2}:\on{Rep}_c(\wh\BG)\to \BVect$
is isomorphic to $\on{Hom}_{\on{Rep}_{c'}(\wh\BG')}(M_c(\BG,\BH_1),M_c(\BG,\BH_2))$.

\section{Proof of \propref{symmetry}}      \label{proof of symmetry}

\ssec{}

We will repeatedly use the following construction:

Let $Z_1\to Z_2$ be a map of schemes of finite type over $\bK$, such that $Z_1$
is principal bundle with respect to a smooth unipotent group-scheme $H$ on $Z_2$.
Let ${\mathcal L}$ be the line bundle on $Z_2$, given by $z\mapsto \on{det}(\fh_z)$,
where $\fh_z$ is the fiber at $z\in Z_2$ of the sheaf of Lie algebras corresponding
to $H$. Let $\wh{Z}_1$ be the total space of the pull-back of the resulting $G_m$-torsor
to $Z_1$.

\begin{lem}  \label{integration along fibers}
Under these circumstances we have a natural map
$$\left(\on{Funct}^{lc}_c(\wh{\bZ}_1)\otimes \BC\right)_{\bG_m}\to
\on{Funct}^{lc}_c(\bZ_2),$$
where $\bG_m$ acts on $\BC$ via the standard character
$\bG_m\to \BZ\overset{1\mapsto q}\to \BC^*$.
\end{lem}

\ssec{}

Let us recall the construction of $M_c(\BG)$, following \cite{GK}, Sect. 5.
To simplify the exposition, we will first assume that $c=1$, in which case
we will sometimes write $M(\BG)$ instead of $M_c(\BG)$.

Consider the set of pairs $(i,Y)$, where $Y$ is a sub-scheme of $G((t))$, stable
under the right action of the congruence subgroup $G^i$. Note that in this case the quotient
$Y/G^i$ is a scheme of finite type over $\bK$.

The above set is naturally filtered: $(i,Y)<(i',Y')$ if $i'\geq i$ and $Y\subset Y'$.
Note also that $Y/G^{i'}\to Y/G^i$ is a principal bundle with respect to
the group $G^i/G^{i'}$.

Let $\bY/\bG^i$ denote the object of $\bSet$, corresponding to the scheme
$Y/G^i$. Consider the
vector space $\bV(i,Y):=\on{Funct}^{lc}_c(\bY/\bG^i)\otimes \mu(\bG[[t]]/\bG^i)$, cf. \cite{GK}, Sect. 3.2,
where for a locally compact group $\bH$, we denote by $\mu(\bH)$ the space of
left-invariant Haar measures on it.

Whenever $(i,Y)<(i',Y')$, we have a natural map $\bV(i',Y')\to \bV(i,Y)$. It is defined
as the composition of the restriction map
$\on{Funct}^{lc}_c(\bY'/\bG^{i'})\to \on{Funct}^{lc}_c(\bY/\bG^{i'})$, followed by
the map
$$\on{Funct}^{lc}_c(\bY/\bG^{i'})\otimes \mu(\bG^i/\bG^{i'})\to
\on{Funct}^{lc}_c(\bY/\bG^i),$$
coming from \lemref{integration along fibers}, using 
$\mu(\bG[[t]]/\bG^{i'})\simeq \mu(\bG[[t]]/\bG^i)\otimes \mu(\bG^i/\bG^{i'})$.

\medskip

We have:
$$M(\BG)=\underset{\underset{(i,Y)}\longleftarrow}{"lim"}\, \bV(i,Y),$$
as a pro-vector space.

\medskip

Let us now describe the action of $\BG\times \wh{\BG}_0$ on $M(\BG)$. For
our purposes it would suffice to do so on the level of
groups of $\bK$-valued points of the corresponding group-indschemes.

For ${\mathbf g}\in G((t))(\bK)$ acting on $M(\BG)$ {\it on the left}, we define
$\bV(i,Y)\to \bV(i,{\mathbf g}\cdot Y)$ to be the natural map. In this way we obtain
an action of ${\mathbf g}$ on the entire inverse system.

To define the right action, for $(i,Y)$ as above, let $j$ be a large enough integer, so that
$\on{Ad}_{{\mathbf g}^{-1}}(G^j)\subset G^i$. Then the right multiplication by ${\mathbf g}$
defines a map of schemes,
$$Y/G^j\to Y\cdot{\mathbf g}/G^i,$$
such that the former is a principal $G^i/\on{Ad}_{{\mathbf g}^{-1}}(G^j)$-bundle over the
latter.

A lift of ${\mathbf g}$ to a point
$\wh{\mathbf g}$ of the central extension $\wh{G}_0$ defines an identification
$\mu(\bG[[t]]/\bG^j)\simeq \mu(\bG[[t]]/\on{Ad}_{{\mathbf g}^{-1}}(G^j))$. Hence,
by \lemref{integration along fibers}, we obtain a map
$$\bV(j,Y)\to \bV(i,Y\cdot{\mathbf g}),$$
and, hence, an action of $\wh{\mathbf g}$ on the inverse system.

\ssec{}

Let now $\wh\BG$ and $c$ be general. We modify the above construction as
follows. For each $Y\subset G((t))$ as above, let $\wh{Y}$ be its pre-image in $\wh{G}$.
Set
$$\bV_c(j,Y):=\left(\on{Funct}^{lc}_c(\bG^j\backslash \wh{\bY})\otimes \BC\right)_{\bG_m}
\otimes \mu(\bG[[t]]/\bG^j),$$
where $\bG_m$ acts naturally on $\wh{\bY}$ and by the character $c$ on $\BC$.
We have:
$$M_c(\BG)=\underset{\underset{(j,Y)}\longleftarrow}{"lim"}\, \bV_c(j,Y),$$
and the action of $\wh\BG\times \wh\BG'$ is described
in the same way as above.

\medskip

By definition, the representation $\wt{M}_{c'}(\BG)$ is the same as
$M_{c'}(\BG)$, viewed as a representation of 
$\wh\BG\times \wh\BG'\simeq \wh\BG'\times \wh\BG$.
Explicitly it can be written down as follows. Consider the set of pairs $(j,Y)$, where $Y\subset G((t))$
is stable under the action of $G^{j}$ {\it on the left}; let $\wh{Y}'$ be the preimage of $Y$
in $\wh{G}'$. We have:
$$\wt{M}_{c'}(\BG)=\underset{\underset{(j,Y)}\longleftarrow}{"lim"}\, \wt{\bV}_{c'}(j,Y),$$
where
$$\wt{\bV}_{c'}(j,Y):=\left(\on{Funct}^{lc}_c(\bG^j\backslash \wh{\bY}')\otimes \BC\right)_{\bG_m}
\otimes \mu(\bG[[t]]/\bG^j),$$
where $\bG_m$ acts naturally on $\wh{\bY}$ and by the character $c'$ on $\BC$.
In this presentation, the {\it right} action of $\wh\BG'$ is defined in an evident fashion, and the {\it left}
action of $\wh\BG$ is defined as in the case of the right action of $\wh\BG_0$ on $M(\BG)$.

\ssec{}

We shall now construct the sought-for map $\wt{M}_{c'}(\BG)\to M_c(\BG)$.
Let us mention that when $G$ is the multiplicative group $G_m$ the sought-for
isomorphism amounts to simply to the inversion on the group.

\medskip

For a pair $(i,Y)$ as in the definition of $M_c(\BG)$, there exists an integer $j$
large enough so that $\on{Ad}_{y^{-1}}(G^j)\subset G^i$ for $y\in Y(\ol{\bK})$.
In particular, over $Y/G^i$ we obtain a group-scheme, denoted
$G^{i,j}_Y$, whose fiber over $y\in Y$ is $G^i/\on{Ad}_{y^{-1}}(G^j)$, and we have a map
\begin{equation} \label{map of quotients}
G^j\backslash Y\to Y/G^i,
\end{equation}
such that the former scheme is a principal $G^{i,j}_Y$-bundle over the
latter.

Note that the fiber of $\wh{Y}$ over a given point $y\in Y$ identifies with
$\on{det}(\on{Ad}_{y}(\fg[[t]]),\fg[[t]])$, where $\fg$ is the Lie algebra of $G$.
Hence, we obtain a natural map
$$\wt{\bV}_{c'}(j,Y)\to \bV_c(i,Y)$$
from \lemref{integration along fibers}.

\medskip

Thus, we obtain a map $\wt{M}_{c'}(\BG)\to M_c(\BG)$, and
from the construction, it is clear that this map respects the action
of $\wh{G}(\bK)\times \wh{G}'(\bK)$. Now \lemref{2.13}
implies that the constructed map is a morphism of
$\wh\BG\times \wh{\BG}'$-representations.

The map in the opposite direction: $M_c(\BG)\to \wt{M}_{c'}(\BG)$
is constructed similarly, and by the definition of the transition maps
giving rise to the inverse systems $M_c(\BG)$ and $\wt{M}_{c'}(\BG)$,
it is clear that both compositions
$M_c(\BG)\to \wt{M}_{c'}(\BG)\to M_c(\BG)$ and
$\wt{M}_{c'}(\BG)\to M_c(\BG)\to \wt{M}_{c'}(\BG)$
are the identity maps.

This proves point (1) of \propref{symmetry}.

\ssec{}

Following \cite{GK}, let us denote by $M(\bG[[t]])$ the pro-vector space
$$\underset{\longleftarrow}{"lim"}\, \on{Funct}_c^{lc}(\bG[[t]]/\bG^i)\otimes \mu(\bG[[t]]/\bG^i),$$
where the transition maps are given by fiber-wise integration. This space
carries an action of the group $\bG[[t]]\times \bG[[t]]$. The convolution product defines
an isomorphism
\begin{equation} \label{conv}
(M(\bG[[t]])\otimes M(\bG[[t]])_{\bG[[t]]}\simeq M(\bG[[t]]),
\end{equation}
where $\bG[[t]]$ acts diagonally.

By construction, as a representation of $\wh\BG$ under the left action,
$M_c(\BG)$ identifies with $i^{\wh\BG}_{\bG[[t]]}(M(\bG[[t]]))$. Therefore,
\begin{equation} \label{map one}
\Hom_{\on{Rep}_c(\wh\BG)}(\wt{M}_{c'}(\BG),M_c(\BG))\simeq
\Hom_{\bG[[t]]}(\wt{M}_{c'}(\BG),M(\bG[[t]])).
\end{equation}
The map $\wt{M}_{c'}(\BG)\to M_c(\BG)$
constructed above corresponds to the natural restriction morphism
$\wt{M}_{c'}(\BG)\to M(\bG[[t]])$.

\medskip

\noindent{\it Remark.}
From the latter description it is not immediately clear why this 
map is compatible with the right $\wh{\BG}'$-action.

Note also that the map $\wt{M}_c(\BG)\to M_{c'}(\BG)$ can be described by
a similar adjunction property with respect to the right $\wh\BG'$-action.

\medskip

Let us prove now point (2) of \propref{symmetry}. For any $\Pi$, which is a representation
of $\wh{\BG}\times \wh{\BG}'$ at levels $(c,c')$ we have:

$$\Hom_{\wh\BG\times \wh{\BG}'}\left(\Pi, (M_c(\BG)\otimes \wt{M}_{c'}(\BG))_{\bG[[t]]}\right)
\simeq \Hom_{\bG[[t]]\times \bG[[t]]}(\Pi,M(\bG[[t]]),$$
with the isomorphism being given by the restriction map
$$(M_c(\BG)\otimes \wt{M}_{c'}(\BG))_{\bG[[t]]}\to (M(\bG[[t]])\otimes M(\bG[[t]]))_{\bG[[t]]},$$
followed by the map of \eqref{conv}.

Let us apply this to $\Pi=M_c(\BG)$. It is clear that both maps appearing in \propref{symmetry}(2),
correspond under the above isomorphism to the restriction map
$M_c(\BG)\to M(\bG[[t]])$. Therefore, these two maps coincide.

\section{Distributions on a stack}   \label{stacks}

\ssec{}

First, let $\bX$ be a locally compact object of $\bSet$. Recall that
${\mathbb Funct}^{lc}(\bX)$ denotes the corresponding (strict) object in
$\BVect$ (cf. \cite{GK}, Sect. 3.2), and $\on{Funct}^{lc}(\bX)=limProj\,
{\mathbb Funct}^{lc}(\bX)$.
The vector space $\on{Distr}_c(\bX)$ introduced in \secref{various distr}
identifies with $\on{Hom}_{\BVect}({\mathbb Funct}^{lc}(\bX),\BC)$,
or, which is the same, with the space of
linear functionals $\on{Funct}^{lc}(\bX)\to \BC$, continuous
in the topology of projective limit.

\medskip

Suppose now that $\bX=X(\bK)$, where $X$ is a {\it smooth} algebraic variety
over $\bK$. In this case we can introduce the subspace $\on{Distr}^{lc}_c(\bX)$
of locally constant distributions on $\bX$ (cf., e.g.,  \cite{GK}, Sect. 5.1).

Indeed, it is well-known that a choice of a top differential form $\omega$
on $X$ defines a measure $\mu(\omega)$ on $\bX$, i.e., a functional on the space
$\on{Funct}^{lc}_c(\bX)$. For $\omega'=\omega\cdot f$, where $f$ is an invertible
function on $X$, we have: $\mu(\omega')=\mu(\omega)\cdot |f|$. Hence, the
subset of elements in $\on{Distr}_c(\bX)$, which can be (locally) written
as $\mu(\omega)\cdot g$, where $g$ is a locally constant function on $\bX$
with compact support, is independent of the choice of $\omega$. This
subset is by definition $\on{Distr}^{lc}_c(\bX)$.

Although the following is well-known, we give a proof for the
sake of completeness:

\begin{prop}  \label{stability of loc con distr}
Let $f:X_1\to X_2$ be a smooth map between smooth varieties over $\bK$. Then

\smallskip

\noindent{\em (1)}
The push-forward map $\on{Distr}_c(\bX_1)\to \on{Distr}_c(\bX_2)$
sends $\on{Distr}^{lc}_c(\bX_1)$ to $\on{Distr}^{lc}_c(\bX_2)$.

\smallskip

\noindent{\em (2)}
If $X_1(\bK)\to X_2(\bK)$ is surjective, then
$f_!:\on{Distr}^{lc}_c(\bX_1)\to\on{Distr}^{lc}_c(\bX_2)$ is
also surjective.

\end{prop}

\begin{proof}

Statement (1) is local in the analytic, and a fortiori in the Zariski topology
on $\bX_1$. Therefore, we can assume that our morphism $f$ factors as
$X_1\overset{f'}\to X_2\times Z\overset{f''}\to X_2$, where $Z$ is another
smooth variety, with $f'$ being \'etale, and $f''$ being the projection on the
first factor.

Since an \'etale map induces a local isomorphism in the analytic
topology,  it is clear that $f'_!$ maps
$\on{Distr}^{lc}_c(\bX_1)$ to $\on{Distr}^{lc}_c(\bX_2\times \bZ)$.

From the definition of $\on{Distr}^{lc}_c(\cdot)$, it is clear that
\begin{equation} \label{distr on prod}
\CD
\on{Distr}^{lc}_c(\bZ_1)\otimes \on{Distr}^{lc}_c(\bZ_2)  @>{\sim}>>
\on{Distr}^{lc}_c(\bZ_1\times \bZ_2)   \\
@VVV  @VVV \\
\on{Distr}_c(\bZ_1)\otimes \on{Distr}_c(\bZ_2) @>>>
\on{Distr}_c(\bZ_1\times \bZ_2).
\endCD
\end{equation}

So the map $f''_!:\on{Distr}^{lc}_c(\bX_2\times \bZ)\to \on{Distr}_c(\bX_2)$
can be identified with
$$\on{Distr}^{lc}_c(\bZ)\otimes \on{Distr}^{lc}_c(\bX_2) \overset{\int\times \on{id}}\to
\on{Distr}^{lc}_c(\bX_2),$$
implying assertion (1) of the proposition.

\medskip

We will prove a slight strengthening of assertion (2). Note that since $f$
is smooth, the image of $\bX_1$ in $\bX_2$ is open, and hence, also closed
in the analytic topology. We will show that $f_1$ maps $\on{Distr}^{lc}_c(\bX_1)$
surjectively onto the subspace of $\on{Distr}^{lc}_c(\bX_2)$, consisting of
distributions, supported on the image.

The assertion is local in the analytic topology on $\bX_2$. Let
$x_2\in X_2(\bK)$ be a point, and let $x_1\in X_1(\bK)$ be some its
pre-image. Then the local factorization of $f$ as $f''\circ f'$ as
above makes the assertion manifest.

\end{proof}

\ssec{}

In what follows we will need a relative version of the above notions. For a smooth morphism
$g:X\to Z$ let $\omega_{rel}$ be a relative top differential form
on $X$. It defines a relative measure $\mu(\omega_{rel}):\on{Funct}^{lc}_c(\bX)\to
\on{Funct}^{lc}_c(\bZ)$. As in the absolute situation, by multiplying $\mu(\omega_{rel})$
by locally constant functions on $\bX$, whose support is proper over $\bZ$, we obtain
a pro-vector sub-space inside $\on{Hom}_{\on{Funct}^{lc}(\bZ)}
(\on{Funct}^{lc}_c(\bX),\on{Funct}^{lc}_c(\bZ))$, which we will denote by
${\mathbb Distr}^{lc}_c(\bX/\bZ)$. Note that when $X=X'\times Z$,
we have: ${\mathbb Distr}^{lc}_c(\bX/\bZ)\simeq \on{Distr}_c^{lc}(\bX')\otimes
{\mathbb Funct}^{lc}(\bZ)$ (the tensor product being taken in the sense of $\BVect$).
We will denote by $\on{Distr}^{lc}_c(\bX/\bZ)$ the vector space
$limProj\, {\mathbb Distr}^{lc}_c(\bX/\bZ)$.

When $f:X_1\to X_2$ is a smooth map of schemes smooth over $Z$, as in \propref{stability of loc con distr}
we have a push-forward map $f_!:{\mathbb Distr}^{lc}_c(\bX_1/\bZ)\to {\mathbb Distr}^{lc}_c(\bX_2/\bZ)$,
which is surjective if $f:X_1(\bK)\to X_2(\bK)$ is; moreover, in this case the map
$f_!:\on{Distr}^{lc}_c(\bX_1/\bZ)\to \on{Distr}^{lc}_c(\bX_2/\bZ)$ is also easily seen to be
surjective. In the particular case when $X_2=Z$ we obtain a map
$\int:{\mathbb Distr}^{lc}_c(\bX/\bZ)\to {\mathbb Funct}^{lc}(\bZ)$.

If $Y$ is another scheme over $Z$, consider the Cartesian diagram
\begin{equation} \label{cart diag}
\CD
X\underset{Z}\times Y @>{f'}>> X \\
@V{g'}VV @V{g}VV  \\
Y @>{f}>> Z.
\endCD
\end{equation}
We have a pull-back map
$f^*:{\mathbb Distr}^{lc}_c(\bX/\bZ)\to {\mathbb Distr}^{lc}_c(\bX\underset{\bZ}\times \bY/\bY)$.

\medskip

Suppose now that the scheme $Z$ is itself smooth, and $X$ is smooth over $Z$ as above.
In this case the spaces
$\on{Distr}^{lc}_c(\bX)$ and $\on{Distr}^{lc}_c(\bZ)$ are well-defined, and we have
an isomorphism
$$\on{Distr}^{lc}_c(\bX)\simeq \on{Distr}^{lc}_c(\bX/\bZ)\underset{\on{Funct}^{lc}(\bZ)}
\otimes \on{Distr}^{lc}_c(\bZ).$$
If in the situation of \eqref{cart diag} $Y$ is also smooth over $\bZ$, and $\xi_Y\in \on{Distr}_c^{lc}(\bY)$,
$\xi_{X/Z}\in \on{Distr}_c^{lc}(\bX/\bZ)$, consider the element $f^*(\xi_{X/Z})\otimes \xi_Y\in
\on{Distr}^{lc}_c(\bX\underset{\bZ}\times \bY)$. We have:
\begin{align} \label{base change for distr}
&f'_!(f^*(\xi_{X/Z})\otimes \xi_Y)=\xi_{X/Z}\otimes f_!(\xi_Y)\in  \on{Distr}_c^{lc}(\bX), \text{ and }\\
&g'_!(f^*(\xi_{X/Z})\otimes \xi_Y)=f^*(g_!(\xi_{X/Z}))\cdot \xi_Y\in \on{Distr}_c^{lc}(\bY).
\end{align}

Finally, let us assume that both maps $f$ and $g$
induce surjections on the level of $\bK$-valued points.

\begin{lem}  \label{Cart product}
The maps $f_!,g_!$ induce an isomorphism
$$\on{Distr}^{lc}_c(\bZ)\simeq \on{coker}\left(\on{Distr}^{lc}_c(\bX\underset{\bZ}\times \bY)
\overset{(f'_!,-g'_!)}\longrightarrow
\on{Distr}^{lc}_c(\bX)\oplus \on{Distr}^{lc}_c(\bY)\right).$$
\end{lem}

\begin{proof}

Let $(\xi_X,\xi_Y)\in \on{Distr}^{lc}_c(\bX)\oplus \on{Distr}^{lc}_c(\bY)$ be an
element such that $f_!(\xi_X)=g_!(\xi_Y)$. We need to find an element
$\xi'\in \on{Distr}^{lc}_c(\bX_1\underset{\bZ}\times \bY)$, such that
$f'_!(\xi')=\xi_X$ and $g_!(\xi')=\xi_Y$. Using \lemref{stability of loc con distr}, we can
assume that $\xi_X=0$.

Let $\xi_{X/Z}$ be an element in $\on{Distr}^{lc}_c(\bX/\bZ)$, such that $\int \xi=1\in \on{Funct}(\bZ)$.
Then $\xi':=f^*(\xi_{X/Z})\otimes \xi_Y$ satisfies our requirements, by \eqref{base change for distr}.

\end{proof}

Let now $X$ and $Y$ be smooth varieties, and $f:Z\times X\to Y$ a map, such that the
corresponding map $f':Z\times X\to Z\times Y$ is smooth.

\begin{lemconstr}  \label{action on distr}
Under the above circumstances we have a natural action map
$$\bZ\times \on{Distr}^{lc}_c(\bX)\to \on{Distr}^{lc}_c(\bY).$$
\end{lemconstr}

\begin{proof}

Consider the map
$$f'_!: {\mathbb Distr}^{lc}_c(\bZ\times \bX/\bZ)\to  {\mathbb Distr}^{lc}_c(\bZ\times \bY/\bZ).$$
By composing it with
$\cdot \otimes 1: \on{Distr}^{lc}_c(\bX)\to {\mathbb Distr}^{lc}_c(\bZ\times \bX/\bZ)$ we obtain a
map
$$\on{Distr}^{lc}_c(\bX)\to \on{Distr}^{lc}_c(\bY)\otimes {\mathbb Funct}^{lc}(\bZ).$$
The latter is, by definition, the same as an action map
$\bZ\times \on{Distr}^{lc}_c(\bX)\to \on{Distr}^{lc}_c(\bY)$.

\end{proof}

\ssec{}

Let $\CY$ be an algebraic stack. We will say that
$\CY$ is $\bK$-admissible (or just admissible) if there exists a smooth
covering $Z\to \CY$, such that for any map $X\to \CY$, the corresponding
map of schemes $$X\underset{\CY}\times Z\to X$$
is surjective on the level of $\bK$-points.

If $\CY$ is admissible, a covering $Z\to \CY$ having the
above property will be called admissible.
It is clear that the class of admissible coverings is closed under
Cartesian products. It is also clear that if $\CY$ is admissible, and
$\CY'\to \CY$ is a representable map, then $\CY'$ is also admissible.

\begin{lem}  \label{admissible stacks}
Suppose that $\CY$ is a stack, which is locally in the Zariski
topology has the form $Z/G$, where $Z$ is a scheme,
and $G$ is an affine algebraic group. Then $\CY$ is admissible.
\end{lem}

\begin{proof}

First, we can assume that $G=GL_n$. Indeed, by assumption,
there is an embedding $G\to GL_n$, and consider the scheme
$Z':=Z\underset{G}\times GL_n$. Then $\CY=Z'/GL_n$.

Now the assertion follows from Hilbert's 90: for $y\in \CY(\bK)$ its
pre-image in $Z$ is a $GL_n$-torsor, which is necessarily trivial.

\end{proof}

\medskip

From now on, we will assume that $\CY$ is admissible. Assume in addition
that $\CY$ is smooth.
We will now define the space, denoted, $\on{Distr}^{lc}_c(\bY)$, of locally constant
compactly supported distributions on $\bY$.

Namely, given two admissible coverings $Z_1,Z_2\to \CY$ we define
$$\on{Distr}^{lc}_c(\bY):=
\on{coker}\left(\on{Distr}^{lc}_c(\bZ_1\underset{\bY}\times \bZ_2)\to
\on{Distr}^{lc}_c(\bZ_1)\oplus \on{Distr}^{lc}_c(\bZ_2)\right).$$

\lemref{Cart product}, combined with \propref{stability of loc con distr}(2), implies
that $\on{Distr}^{lc}_c(\bY)$ is well-defined, i.e., is independent of the choice of
$Z_1$, $Z_2$.

If $f:\CY_1\to \CY_2$ is a smooth representable map of (smooth admissible)
stacks, from \propref{stability of loc con distr}(1) we obtain that there
exists a well-defined map $f_!:\on{Distr}^{lc}_c(\bY_1)\to \on{Distr}^{lc}_c(\bY_2)$.

\medskip

Assume now that $\CY=Z/G$, where $G$ is an algebraic group acting on $Z$.
By \lemref{action on distr}, we have an action of $\bG$ on the vector space
$\on{Distr}^{lc}_c(\bZ)$. From \lemref{Cart product} we obtain:

\begin{cor} \label{distr as coinv}
For $\CY$ as above,
$$\on{Distr}^{lc}_c(\bY)\simeq \on{Coinv}_\bG(\on{Distr}^{lc}_c(\bZ)).$$
\end{cor}

\ssec{Relative version}  \label{relative stack version}

Assume now that $\CY$ is a stack, endowed with a smooth map to
a scheme $Z$. For a pair of admissible coverings $X_1,X_2\to \CY$,
we define the pro-vector space ${\mathbb Distr}^{lc}_c(\CY/Z)$ as
$$\on{coker}\left({\mathbb Distr}^{lc}_c(\bX_1\underset{\bY}\times \bX_2/\bZ)\to
{\mathbb Distr}^{lc}_c(\bX_1/\bZ)\oplus {\mathbb Distr}^{lc}_c(\bX_2/\bZ)\right).$$
A relative version of \lemref{Cart product} shows that this is well-defined,
i.e., independent of the choice of $X_1$ and $X_2$.

Finally, the assertion of \lemconstrref{action on distr} remains valid, where
$Z$ is a scheme, $\CY,\CY'$ are smooth stacks, and the map
$f:Z\times \CY\to \CY'$ is such that the corresponding map
$f':Z\times \CY\to Z\times \CY'$ is smooth and representable.

\section{Induction via the moduli stack of bundles}  \label{sect thick}

\ssec{}

Let $X$ be a (smooth complete) algebraic curve over $\bK$,
$\bx\in \bX$ a rational point, and let $t$ be a coordinate near $\bx$.

If $G$ be a split reductive group, let $\Bun_G$ denote the moduli
stack of principal $G$-bundles on $X$. For $i\in \BZ$, let $\Bun^{i,\bx}_G$ 
denote the stack classifying bundles equipped with a trivialization 
on the $i$-th infinitesimal neighbourhood of $\bx$. By construction, 
$\Bun^{i,\bx}_G$ is a principal $G[[t]]/G^i$-bundle over $\Bun_G$.

If $\CY\subset \Bun_G$ is an open sub-stack of finite type, we let
$\CY^{i,\bx}$ denote its pre-image in $\Bun^{i,\bx}_G$. The
following is well-known:

\begin{lem}  \label{stabilization}
For any $\CY\subset \Bun_G$ of finite type and $i$ large enough,
the stack $\CY^{i,\bx}$ is a scheme of finite type.
\end{lem}

For $\CY$ as above, we let $\CY^{\infty,\bx}$ denote the object of
$\on{Pro}(Sch^{ft})$ equal to $\underset{\longleftarrow}{"lim"}\, \CY^{i,\bx}$.
We let $\Bun^{\infty,\bx}_G$ denote the object
$$\underset{\CY}{\underset{\longrightarrow}{"lim"}}\,\CY^{\infty,\bx}\in
\on{Ind}(\on{Pro}(Sch^{ft})).$$

Another basic fact is that $G((t))$, viewed as a group-object of
$\on{Ind}(\on{Pro}(Sch^{ft}))$, acts on $\Bun^{\infty,\bx}_G$ in the sense
of the tensor structure on $\on{Ind}(\on{Pro}(Sch^{ft}))$.

\ssec{}

By \lemref{admissible stacks}, the stacks
$\CY^{i,\bx}$ are admissible. Set $\bW^i_\CY=\on{Distr}^{lc}_c(\bY^{i,\bx})$.
For $\CY_1\hookrightarrow \CY_2$ we have a natural push-forward
map on the level of distributions $\bW^i_{\CY_1}\to \bW^i_{\CY_2}$.
Set
$$\bW^i:=\underset{\CY}{\underset{\longrightarrow}{lim}}\, \bW^i_\CY\in Vect.$$

For a fixed $\CY$ and  $j>i$ we have a smooth representable map of stacks
$\bY^{j,\bx}\to \bY^{i,\bx}$; hence we obtain a map
$\bW^j_\CY\to \bW^i_\CY$ and, finally, a map $\bW^j\to \bW^i$.

We define the pro-vector space
$$\BW_{X,\bx}:=\underset{i}{\underset{\longleftarrow}{"lim"}}\, \bW^i.$$

Now we are ready to state:

\begin{thm}  \label{thick induction}
The pro-vector space $\BW_{X,\bx}$ carries a natural action of the group
$\BG$, such that $\on{Coinv}_{\bG^i}(\BW_{X,\bx})\simeq \bW^i$.
\end{thm}

Note that by construction we have:
\begin{cor}   \label{adm}
The $\BG$-representation $\BW_{X,\bx}$ is admissible.
\end{cor}

Indeed, the coinvariants $\on{Coinv}_{\bG^i}(\BW_{X,\bx})\simeq \bW^i$ all
belong to $Vect$.

\ssec{Proof of \thmref{thick induction}}

Let $G((t))=\underset{k}{\underset{\longrightarrow}{"lim"}}\, Z_k$ with
$Z_k=\underset{l}{\underset{\longleftarrow}{"lim"}}\, Z^l_k$, where
$Z^l_k$ are schemes of finite type.

To define an action
$$\BG\times \BW_{X,\bx}\to \BW_{X,\bx}$$
we need to give for every $k$ and $i$ a map
$$\bZ^l_k\times \bW^j\to \bW^i$$
defined for $j$ and $l$ sufficiently large.

For $k$ and $i$ as above let $j$ be such that
$\on{Ad}_{Z_k}(G^j)\subset G^i$. The action of $G((t))$ on
$\Bun^{\infty,\bx}_G$ yields a map of stacks $Z_k\times
\Bun^{j,\bx}_G\to \Bun^{i,\bx}_G$, which factors through $Z^l_k$ for
some $l$. Moreover, for every sub-stack $\CY\subset \Bun_G$ of
finite type, there exists another sub-stack $\CY'$ of finite type,
such that we have a map
$$Z^l_k\times \CY^{j,\bx}\to \CY'{}^{j,\bx}.$$

We claim that for $i,j,k,l,\CY,\CY'$ as above, we have a map
\begin{equation} \label{required map}
\bZ^l_k\times \bW^j_\CY\to \bW^j_{\CY'}.
\end{equation}
This follows from the stack-theoretic version of
\lemconstrref{action on distr}, cf. \secref{relative stack version}.
The fact that the resulting action map $\BG\times \BW_{X,\bx}\to
\BW_{X,\bx}$ respects the group law on $\BG$ is a straightforward
verification.

\medskip

To compute $\on{Coinv}_{\bG^i}(\BW_{X,\bx})$ note that
$G[[t]]$, and hence all $G^i$, act on each $\Bun^{j,\bx}_G$
individually.

Hence,
$$\on{Coinv}_{\bG^i}(\BW_{X,\bx})\simeq
\underset{j\geq i}{\underset{\longleftarrow}{"lim"}}\,
\on{Coinv}_{\bG^i/\bG^j}(\bW^j).$$

We claim that for $j\geq i$, 
$\on{Coinv}_{\bG^i/\bG^j}(\bW^j)\simeq \bW^i$. Indeed, since
each $\CY^{j,\bx}$ is stable under $G^i/G^j$, we have:
$$\on{Coinv}_{\bG^i/\bG^j}(\bW^j)\simeq
\underset{\CY}{\underset{\longrightarrow}{lim}}\, \on{Coinv}_{\bG^i/\bG^j}(\bW^j_\CY)\simeq
\underset{\CY}{\underset{\longrightarrow}{lim}}\, \bW^i_\CY\simeq \bW^i,$$
where the middle isomorphism follows from \corref{distr as coinv}.

\ssec{Variants and Generalizations}

Recall that the stack $\Bun_G$ is endowed with a canonical line bundle
$\CL_{\Bun_G}$, with the basic property that the $\BG$-action on
$\Bun_G^{\bx,\infty}$ extends to an action of a central extension
$\wh\BG$ on the pull-back of $\CL_{\Bun_G}$ to $\Bun_G^{\bx,\infty}$.

By the same token, we consider now a representation
$\wh\BW_{X,\bx}$ of $\wh\BG$, and for every $c:\bG_m\to \BC^*$
the object
$$\BW_{X,\bx,c}:=(\wh\BW_{X,\bx}\otimes \BC)_{\bG_m}\in \on{Rep}_c(\wh\BG).$$

\medskip

Note that instead of a single point $\bx$ we could have considered
any finite collection $\ol{\bx}=\bx_1,...,\bx_n$ of rational points.
By repeating the construction we obtain a pro-vector space
$\BW_{X,\ol{\bx}}$, acted on by the product $\underset{k}\Pi\,
\wh{\BG}_{\bx_k}$, where each $\wh{\BG}_{\bx_k}$ identifies with
$\wh{\BG}$ once we identify the local ring of $\bX$ at $\bx_k$ with
${\mathbf F}$.

Again, for a choice of a character $c:\bG_m\to \BC^*$, we obtain a
representation of $\underset{k}\Pi\, \wh{\BG}_{\bx_k}$, denoted $\BW_{X,\ol{\bx},c}$,
such that the center $\bG^k_m$ acts via the multiplication map
$\bG^k_m\to \bG_m$.

\ssec{}

From now on we will suppose that $X$ is isomorphic to the projective
line $P^1$, and the number of points is two, which we will denote by
$\bx_1$, and $\bx_2$, respectively. Assume also that $G$ is semi-simple and
simply connected.

Consider the representation
$\wh\BW_{P^1,\bx_1,\bx_2,c}$ of $\wh\BG_{\bx_1}\times \wh\BG_{\bx_2}$.
Let us take its coinvariants with respect to 
$\bI^{00}_{\bx_1}\subset \wh\BG_{\bx_1}$.
By \corref{action on coinv}, on the resulting pro-vector space we
will have an action of $\wh\BG_{\bx_2}$; we will denote this representation
by $\Pi^{\on{thick}}_c$, i.e.,
$$\Pi^{\on{thick}}_c=\on{Coinv}_{\bI^{00}_{\bx_1}}(\wh\BW_{P^1,\bx_1,\bx_2,c}).$$

By \thmref{endomorphisms}, the algebra $\overset{\cdot\cdot}\sH_{q,c'}$ acts
on $\Pi^{\on{thick}}_c$ by endomorphisms. Consider now
$$\bU_c:=\on{Coinv}_{\bI^{00}_{\bx_1}\times 
\bI^{00}_{\bx_2}}(\wh\BW_{P^1,\bx_1,\bx_2,c})\simeq
\on{Coinv}_{\bI^{00}_{\bx_2}}(\Pi^{\on{thick}}_c).$$

By \corref{adm}, this is a vector space, endowed with two commuting
actions of $\overset{\cdot\cdot}\sH_{q,c'}$. We have:

\begin{thm}   \label{regular}
There exists a canonically defined vector $\one_{\bU_c}\in \bU_c$, which
freely generates $\bU_c$ under each of the two
$\overset{\cdot\cdot}\sH_{q,c'}$-actions.
\end{thm}

\section{Proof of \thmref{regular}}  \label{sect regular}

\ssec{}

Let $W_{aff}$ be the affine Weyl group corresponding to $G$.
Since $G$ was assumed simply connected, $W_{aff}$ is a Coxeter
group.

If $\alpha$ is a simple affine root, let $I_\alpha\subset \wh{G}$
denote the corresponding sub-minimal parahoric; let $N(I_\alpha)$
denote the (pro)-unipotent radical of $I_\alpha$, and
$M_\alpha:=I_\alpha/N(I_\alpha)$ the Levi quotient.

By definition, $M_\alpha$ is a reductive group of semi-simple rank
$1$, with a distinguished copy of $G_m$ in its center; we will
denote by $M_\alpha'$ the quotient $M_\alpha /G_m$.
Let $B_\alpha$ denote the Borel subgroup of $M_\alpha$, and
$\bB^0_\alpha$ the kernel of
$$\bB_\alpha\to \bT\to \Lambda.$$
Let $\Pi^\alpha_c$ be the quotient of the principal series representation of
$\bM_\alpha$, given by the condition that $\bG_m\subset \bM_\alpha$
acts by the character $c$, i.e.,
\begin{equation} \label{Pi alpha}
\Pi^\alpha_c=\left(\on{Funct}^{lc}_c(\bM_\alpha/\bB^0_\alpha)\otimes \BC\right)_{\bG_m}.
\end{equation}

Let us denote by
$\overset{\cdot}\sH{}^\alpha_{q,c}$ the corresponding affine Hecke algebra
of $M_\alpha$, i.e., the algebra of endomorphisms of the functor
$\on{Coinv}_{\bB^0_\alpha}:\on{Rep}(\bM_\alpha,Vect)_c\to Vect$,
or which is the same, the algebra of endomorphisms of $\Pi^\alpha_c$
as a $\bM_\alpha$-representation. It is well-known that
$\bU_\alpha:=\on{Coinv}_{\bB^0_\alpha}(\Pi^\alpha_c)$, as a bi-module
over $\overset{\cdot}\sH{}^\alpha_{q,c}$, is isomorphic to the regular
representation. In a sense, \thmref{regular} generalizes this result to the affine
case.

\medskip

The functor $\Pi\mapsto \on{Coinv}_{\bI^{00}}(\Pi)$ on
$\on{Rep}_c(\wh\BG)$ can be factored into two steps. We
first apply the functor
$$r^{\wh{\BG}}_{\bM^\alpha}:\on{Rep}_c(\wh\BG)\to \on{Rep}_c(\bI_\alpha,\BVect)
\overset{\on{Coinv}_{\bN_\alpha}}
\longrightarrow \on{Rep}_c(\bM_\alpha,\BVect),$$
where the first arrow is the forgetful functor, and then apply
$$\on{Coinv}_{\bB^{0}_\alpha}:
\on{Rep}_c(\bM_\alpha,\BVect)\to \BVect.$$
In particular, endomorphisms of the latter functor map
to endomorphisms of the composition. As a result, we
obtain the canonical embedding $\overset{\cdot}\sH{}^\alpha_{q,c'}\to
\overset{\cdot\cdot}\sH_{q,c}$.

Recall also that the group-algebra $\BC[\Lambda]$ is canonically a subalgebra
in $\overset{\cdot\cdot}\sH_{q,c}$, 
contained in each $\overset{\cdot}\sH{}^\alpha_{q,c'}$.

\ssec{}

The strategy of the proof of \thmref{regular}
will be as follows. We will endow the vector space
$\bU_c$ with an increasing filtration
$$\bU_c=\underset{w\in W_{aff}}\cup\, \bU_w$$ with
$\bU_{w_1}\subset \bU_{w_2}$ if and only if $w_1\leq w_2$ in
the Bruhat order. This filtration will be stable under the action
of $\BC[\Lambda]\subset \overset{\cdot\cdot}\sH_{q,c}$ with
respect to both actions of the latter on $\bU_c$.

The subquotients
$$\bU^w:=\bU_w/\underset{w'<w}\cup\, \bU_{w'}$$
will be free $\Lambda$-modules of rank $1$ (with respect to
each of the actions of $\overset{\cdot\cdot}\sH_{q,c}$). In particular,
for $w=1$, the space $\bU^1\simeq \bU_1$ will contain a canonical
element $\one_{\bU^1}\in \bU^1$, which generates $\bU^1$ under
each of the $\Lambda$-actions. This will be the element
$\one_{\bU_c}$ of \thmref{regular}.

\medskip

Moreover, the following crucial property will be satisfied.
Suppose that $w$ is an element of $W_{aff}$, and $s_\alpha$
is a simple affine reflection, such that $s_\alpha\cdot w>w$
(resp., $w\cdot s_\alpha>w$). Then the subquotient
\begin{equation} \label{subquotient}
\bU_{s_\alpha\cdot w}/\underset{w'<s_\alpha\cdot w,w'\neq w}\cup\, \bU_{w'}, \text{ (resp., }
\bU_{w\cdot s_\alpha}/\underset{w'<w\cdot s_\alpha,w'\neq w}\cup\, \bU_{w'})
\end{equation}
is stable under the action of $\overset{\cdot}\sH{}^\alpha_{q,c}$, embedded
into the first (resp., second) copy of $\overset{\cdot\cdot}\sH_{q,c}$,
and as a $\overset{\cdot}\sH{}^\alpha_{q,c}$-module, it is isomorphic
to $\overset{\cdot}\sH{}^\alpha_{q,c}\underset{\BC[\Lambda]}\otimes \bU^w$.

The existence of a filtration with the above properties clearly
implies the assertion of the theorem.

\ssec{}

Let $I^0$ denote the (pro)-unipotent radical of $I$; we have
$\bI^{00}/\bI^0\simeq \bT^0$, where $\bT^0\subset \bT$ is the maximal
compact subgroup of $\bT$.

Consider the scheme $\CG_G:=\Bun_G^{\infty,\bx_1,\bx_2}/I_{\bx_2}$,
called the thick Grassmannian of $G$. By definition, it classifies
principal $G$-bundles on $X=P^1$, endowed with a trivialization at
the formal neighbourhood of $\bx_1$ and a reduction to $B$ of their
fiber at $\bx_2$. Consider also the base affine space
$\wt{\CG}_G:=\Bun_G^{\infty,\bx_1,\bx_2}/I^0_{\bx_2}$, which is a
principal $T$-bundle over $\CG_G$. The loop group $G((t))$, where
$t$ is the coordinate near $\bx_1$ acts naturally on both $\CG_G$
and $\wt{\CG}_G$.

It is well-known that $\CG_G$ can be written as a union of open
sub-schemes $\CG_{G,w}$, $w\in W_{aff}$, each being stable under the
action of $I_{\bx_1}=I\subset G((t))$, such that $\CG_{G,w_1}\subset
\CG_{G,w_2}$ if and only if $w_1<w_2$ in the Bruhat order. Let us
denote by $\CG^w_G$ the locally closed sub-scheme
$\CG_{G,w}-\underset{w'<w}\cup\,  \CG_{G,w'}$, and by
$\wt\CG_{G,w}$, $\wt\CG^w_G$ the corresponding sub-schemes in
$\wt\CG$.  It is well-known that the group $I^0$ (resp., $I$) acts
transitively on each $\CG^w_G$ (resp., $\wt\CG^w_G$) with
finite-dimensional unipotent stabilizers. Choosing a point in each
$\wt\CG^w_G$, we will denote by $N_w$ its stabilizer in $I$, or,
which is the same, the stabilizer in $I^0$ of the projection of this point to
$\CG^w_G$.

\medskip

Consider the stack
$$\Bun_G^{\bx_1,\bx_2}:=
\Bun_G^{\infty,\bx_1,\bx_2}/(I^0_{\bx_1}\times I^0_{\bx_2})\simeq
\wt\CG_G/I^0.$$ By definition, it classifies $G$-bundles on $X=P^1$
with a reduction to the maximal unipotent at $\bx_1$ and $\bx_2$,
and it carries a natural action of the group $T\times T$. From the
above discussion, we obtain that $\Bun_G^{\bx_1,\bx_2}$ can be
canonically written as a union of open sub-stacks of finite type
$$\Bun_G^{\bx_1,\bx_2}=\underset{w\in W_{aff}}\cup\,\CY_w$$
with $\CY_{w_1}\subset \CY_{w_2}$ if and only if $w_1\leq w_2$.

\medskip

Consider the locally-closed sub-stack
$\CY^w:=\CY_w-\underset{w'<w}\cup\, \CY_{w'}$. We obtain that
$\CY^w$ is isomorphic to $T\times (\on{pt}/N_w)$, where $N_w$ is as
above. The first copy of $T$ acts via multiplication on the first
factor, and the action of the second copy is twisted by the
projection of $w$ to the finite Weyl group, acting by automorphisms
on $T$.

We will denote by $\wh\CY_w$, $\wh\CY^w$ the pull-back of the total
space of the $G_m$-torsor corresponding to $\CL_{\Bun_G}$ to these
sub-stacks.

\ssec{}

We have:
$$\bU_c\simeq \left(\underset{\underset{w}\longrightarrow}{lim}\,
\on{Distr}^{lc}_c(\wh\CY_w)\otimes \BC\right)_{\bT^0\times \bT^0\times \bG_m}.$$
Set
\begin{equation}
\bU_w:=(\on{Distr}^{lc}_c(\wh\CY_w)\otimes \BC)_{\bT^0\times \bT^0\times \bG_m}.
\end{equation}

We claim that each $\bU_w$ maps injectively into $\bU_c$; and the images of
$\bU_w$ define a filtration
with the required properties. One thing is clear, however: by construction,
$\bU_w$ carries an action of $\Lambda\times \Lambda$, and its map to
$\bU_c$ is compatible with this action.

\ssec{}

To proceed we need to introduce some more notation.
Let $Z$ be a smooth scheme, and let $\CL$ be a line bundle on
$Z$. Let $\overset{\circ}\CL$ denote the total space
of the corresponding $G_m$-torsor over $Z$.
We will denote by $\on{Distr}^{lc}_c(Z)_{\CL}$ the space
$$\left(\on{Distr}^{lc}_c(\overset{\circ}\bL)\otimes \BC\right)_{\bG_m},$$
where $\bG_m$ acts on $\BC$ via the standard character
$\bG_m\to \BZ\overset{1\mapsto q}\longrightarrow \BC^*$.

Let now $Z_1\subset Z$ be a smooth closed sub-scheme, and let
$Z_2$ be its complement. We have:

\begin{lem}   \label{SES of distr}
There exists a natural short exact sequence:
$$0\to \on{Distr}^{lc}_c(Z_2)_{\CL}\to
\on{Distr}^{lc}_c(Z)_{\CL}\to \on{Distr}^{lc}_c(Z_1)_{\CL\otimes \CL_n}\to 0,$$
where $\CL_n$ is the top power of the normal bundle
to $Z_1$ inside $Z$.
\end{lem}

\begin{proof}

Note that by definition we have 
$$\on{Distr}^{lc}_c(Z)_{\CL_0}\simeq \on{Funct}^{lc}_c(Z),$$
where $\CL_0$ is the inverse of the line bundle of top forms
on $Z$. The assertion of the lemma follows now from the fact
that for any $Z_1\subset Z$
we have a short exact sequence for the corresponding spaces of
locally constant functions with compact support:
$$0\to \on{Funct}^{lc}_c(Z_2)\to
\on{Funct}^{lc}_c(Z)\to \on{Funct}^{lc}_c(Z_1)\to 0.$$

\end{proof}

For each $w<w'$, the open embedding $\CY_w\hookrightarrow \CY_{w'}$ can
be covered by an open embedding of schemes $Z_w\hookrightarrow Z_{w'}$, such that
$\CY_{w'}=Z_{w'}/N$, $\CY_w=Z_w/N$with $N$ being a unipotent algebraic group.
Therefore, by \lemref{distr as coinv} and the exactness of the functor $\on{Coinv}_\bN$,
the map $\bU_w\to \bU_{w'}$ is an embedding. Hence, $\bU_w\to \bU_c$ is also
an embedding.

Moreover, we claim that from \lemref{SES of distr} we obtain a (non-canonical) isomorphism
\begin{equation} \label{sub-quotient}
\bU^w\simeq(\on{Distr}^{lc}_c(\CY^w)_{\bT^0\times \bT^0},
\end{equation}
compatible with the $\Lambda\times \Lambda$-action. 

\medskip

Indeed, {\it a priori}, $\bU^w\simeq\left((\on{Distr}^{lc}_c(\CY^w)_\CL\right)_{\bT^0\times \bT^0}$
for a certain $T\times T$-equivariant line bundle $\CL$ on $\CY^w$. However, from
the description of $\CY^w$ as $T\times (\on{pt}/N_w)$, this line bundle is (non-canonically)
trivial. Note, however, that this line bundle is canonically trivial for $w=1$.

Now, the same description of $\CY^w$ implies that $\on{Distr}^{lc}_c(\CY^w)\simeq
\on{Funct}^{lc}_c(\bT)$, with the first action of $\bT$ being given by multiplication, and
the second action is twisted by $w$. This implies that 
$\bU^w\simeq (\on{Funct}^{lc}_c(\bT))_{\bT^0\times \bT^0}\simeq \BC[\Lambda]$.

\ssec{}

We will now study the subquotient 
$\bU_{s_\alpha\cdot w}/\underset{w'<s_\alpha\cdot w, w'\neq w}\cup\, \bU_{w'}$,
where $s_\alpha$ is a simple affine reflection such that $s_\alpha\cdot w>w$. (The case
$w\cdot s_\alpha> w$ is analyzed similarly.)

Note first of all that for any $w'\in W_{aff}$, we have $I_\alpha\cdot \CG_G^{w'}\subset
\CG_G^{w'}\cup \CG_G^{s_\alpha\cdot w'}$. Hence, the open subset $\CG_{G,s_\alpha\cdot w}$
is $I_\alpha$-stable, and so is the union $\underset{w'<s_\alpha\cdot w,w'\neq w}\cup\, \CG^{w'}_G$.
Therefore, the subquotient in \eqref{subquotient} is indeed
$\overset{\cdot}\sH{}^\alpha_{q,c'}$-stable.

We will consider two additional stacks. One is
$'\CY:=\Bun_G^{\infty,\bx_1,\bx_2}/(N(I_\alpha)_{\bx_1}\times I^0_{\bx_2})$, on which we
have an action of $M'_\alpha$. We will denote by $pr$ the projection
$$'\CY\overset{pr}\longrightarrow  {}'\CY/N_\alpha\simeq
\Bun_G^{\bx_1,\bx_2},$$
where $N_\alpha:=B_\alpha\cap I^0$.

\medskip

Another stack is the quotient
$$''\CY:=\Bun_G^{\infty,\bx_1,\bx_2}/(M'_\alpha\times I^0).$$

The stack $''\CY$ can be written as a union of open sub-stacks
$''\CY_{w}$ numbered by left cosets $\{1,s_\alpha\}\backslash
W_{aff}$; we will denote by $''\CY^w$ the corresponding locally
closed sub-stacks. Let also $'\CY_w$ and $'\CY^w$ denote the
pre-images of the corresponding sub-stacks in $'\CY$, and
$'\wh{\CY}_w$, $'\wh{\CY}^w$ the total spaces of the $G_m$-torsors,
corresponding to the pull-backs of the line bundle $\CL_{\Bun_G}$.

If $w$ is an element of $W_{aff}$ we have:
\begin{equation} \label{strata}
'\CY^w=pr^{-1}(\CY^w\cup \CY^{s_\alpha \cdot w}).
\end{equation}

\medskip

Using \lemref{SES of distr}, the subquotient \eqref{subquotient}
is isomorphic to
$$\left(\on{Distr}^{lc}_c({}'\wh{\CY}^w)\otimes \BC\right)_{(\bB_\alpha^0)_{\bx_1}\times 
(\bT^0)_{\bx_2}\times \bG_m}.$$
The vector space
$\left(\on{Distr}^{lc}_c({}'\wh{\CY}^w)\otimes \BC\right)_{(\bT^0)_{\bx_2}\times \bG_m}$ is naturally
a representation of the group $\bM_\alpha$. We claim that as such,
\begin{equation} \label{ident of Pi alpha}
\left(\on{Distr}^{lc}_c({}'\wh{\CY}^w)\otimes \BC\right)_{(\bT^0)_{\bx_2}\times \bG_m}\simeq
\Pi^\alpha_c.
\end{equation}

Clearly, the above isomorphism implies our assertion about the action of
$\overset{\cdot}\sH{}^\alpha_{q,c'}$ on the subquotient in \eqref{subquotient}.

\ssec{}

To prove \eqref{ident of Pi alpha} let us observe that $''\CY\simeq {}'\CY/M'_\alpha$
and that $''\CY^w\simeq \on{pt}/N_{w,\alpha}$, where
$N_{w,\alpha}$ is a unipotent group, so that the Cartesian
product
$$\on{pt}\underset{''\CY^w}\times {}'\CY^w$$ is isomorphic to
$M'_\alpha$, and the action of $N_{w,\alpha}$ on $M'_\alpha$
comes from a surjective homomorphism $N_{w,\alpha}\to N_\alpha$
and the action of the latter on $M'_\alpha$ by right multiplication.
Hence,
$$\on{Distr}^{lc}_c({}'\wh\CY^w)\simeq
\left(\on{Distr}^{lc}_c(\bM_\alpha)\right)_{\bN_\alpha}\simeq
\left(\on{Distr}^{lc}_c(\bM_\alpha/\bN_\alpha)\right),$$
implying \eqref{ident of Pi alpha}.

\end{document}